\documentclass[10pt,a4paper,final]{article}
\usepackage[left=2.250cm, right=2.250cm, top=3.750cm, bottom=3.75cm]{geometry}

\usepackage{titlesec}

\usepackage{graphics}
\usepackage{graphicx}
 \usepackage{epsfig}

\usepackage{amssymb}
\usepackage{stmaryrd}
 \usepackage{amsthm}
  \usepackage{amsmath}
\usepackage{tikz-cd}

 \usepackage{float}
 \usepackage{newfloat}
\usepackage{lineno}

\usepackage{ccicons}
\usepackage{epstopdf}

\usepackage{multirow}
\usepackage{color}
\usepackage{hhline}
\usepackage{booktabs}

\usepackage{comment}

\newtheorem{theorem}{Theorem}[section]

\newtheorem{remark}{Remark}[section]

\newtheorem{lemma}{Lemma}[section]

\newtheorem*{weakproblem}{Weak problem}

\definecolor{nverde}{RGB}{0,61,0} 
\definecolor{cr1}{RGB}{200,0,0}
\definecolor{cr2}{RGB}{0,0,200}
\definecolor{cr12}{RGB}{100,0,100}

\newcommand{\xx}{\mathbf{x}}

\newcommand{\ave}[1]{\{\!\!\{#1\}\!\!\}}
\newcommand{\jump}[1]{\llbracket #1 \rrbracket}

\newcommand{\Mach}{\mathrm{M}}

\newcommand{\RT}{\mathrm{RT}}
\newcommand{\DG}{\mathrm{d}\mathcal{P}}

\DeclareMathOperator{\dS}{\mathrm{dS}}

\DeclareMathOperator{\dx}{\mathrm{d\xx}}

\newcommand{\density}{\rho}
\newcommand{\densityh}{\rho_h}
\newcommand{\entropy}{S}
\newcommand{\entropytest}{R}
\newcommand{\entropytesth}{\entropytest_h}
\newcommand{\entropyh}{S_h}
\newcommand{\mom}{\mathbf{m}}
\newcommand{\momh}{\mom_h}
\newcommand{\vel}{\mathbf{u}}
\newcommand{\vor}{\boldsymbol{\omega}}
\newcommand{\pres}{p}
\newcommand{\velh}{\mathbf{u}_h}
\newcommand{\veltesth}{\mathbf{v}_h}
\newcommand{\presh}{p_h}
\newcommand{\prestesth}{q_h}
\newcommand{\vorh}{\boldsymbol{\omega}_h}
\newcommand{\vortest}{\mathbf{z}}
\newcommand{\vortesth}{\mathbf{z}_h}
\newcommand{\Grad}{\nabla}
\newcommand{\Curl}{\nabla\times}
\newcommand{\Div}{\nabla\cdot}

\newcommand{\vis}{\mu}
\newcommand{\tEnd}{t_{\mathrm{end}}}
\newcommand{\ddt}{\partial_t}

\newcommand{\dt}{\Delta t}

\newcommand{\nvec}{\mathbf{n}}

\DeclareFloatingEnvironment[fileext=lop]{Diagram}

\graphicspath{ {./figures/} }

\usepackage{authblk}
\usepackage{fancyhdr}
\usepackage{changepage}
\newcommand{\address}[2]{\affil[#1]{#2}}

\newcommand{\DMUNITN}{1}
\newcommand{\LAM}{2}

\title{An asymptotic-preserving and exactly mass-conservative semi-implicit scheme for weakly compressible flows based on compatible finite elements}

\author[\DMUNITN]{E. Zampa}
\author[\LAM]{M. Dumbser}
\address{\DMUNITN}{Department of Mathematics, University of Trento, Via Sommarive 14, 38123 Trento, Italy}

\address{\LAM}{Laboratory of Applied Mathematics, DICAM, University of Trento, via Mesiano 77, 38123 Trento, Italy} 

\affil[ ]{\footnote{
		\textit{enrico.zampa@unitn.it}  (E. Zampa),
		\textit{michael.dumbser@unitn.it} (M. Dumbser)}}

\date{12 July, 2024}

\newcommand{\keywords}[1]{\textbf{keywords}-- #1}

\allowdisplaybreaks

\begin{document}
\maketitle 

\begin{abstract}
	\textcolor{black}{
We present a novel asymptotic-preserving semi-implicit finite element method for weakly compressible and incompressible flows based on compatible finite element spaces. The momentum is sought in an  $H(\mathrm{div})$-conforming space, ensuring exact pointwise mass conservation at the discrete level. We use an explicit discontinuous Galerkin-based discretization for the convective terms, while treating the pressure and viscous terms implicitly, so that the CFL condition depends only on the fluid velocity. To handle shocks and damp spurious oscillations in the compressible regime, we incorporate an \textit{a posteriori} limiter that employs artificial viscosity and is based on a discrete maximum principle. By using hybridization, the final algorithm requires solving only symmetric positive definite linear systems. As the Mach number approaches zero and the density remains constant, the method converges to an $H(\mathrm{div})$-based discretization of the incompressible Navier-Stokes equations in the vorticity-velocity-pressure formulation. Several numerical tests validate the proposed method.}
\end{abstract}

\keywords{
finite element exterior calculus; compatible finite elements; semi-implicit scheme; weakly compressible flows; incompressible Navier-Stokes equations
}


\section{Introduction}\label{sec:intro}
\color{black}
The numerical approximation of the incompressible Navier-Stokes equations is an active topic of research in which compatible finite elements \cite{Arnold1, Arnold2, Hiptmair} have shown to be natural candidates for high order structure-preserving discretizations, \textcolor{black}{that is, discretizations that conserve  energy, the Hamiltonian structure of the original system, involutions, etc.} Since the pioneering work of Cockburn, Kanschat and Sch\"{o}tzau \cite{CoKaSc2006}, $H(\mathrm{ div})$-based finite element approximations of the incompressible Euler equations \cite{GuShuSe2017, NaCo2018, Gawlik20} have become more and more attractive since they produce exactly divergence-free velocities. Extending these methods to the full Navier-Stokes equations is not a trivial task, since the viscous stress tensor is not bounded in $H(\mathrm{div}$). There are two possible ways to circumvent this obstacle: the first is to use a DG based approximation of the viscous term as in \cite{LeSc16, RhebergenWells18, Fu19}, whereas the second is to introduce additional fields to obtain a conforming term \cite{MEEVC, ZPGR22, Hanot23, CaCPFa2023, MEEVC24}. 

On the other hand, the numerical approximation of compressible flows with compatible finite elements is still at the dawn. In particular, we mention the works by Gawlik and Gay-Balmaz \cite{Gawlik20, Gawlik21, Gawlik24}, the Versatile Mixed Method by Miller and Williams \cite{MillerWilliams2024} and the method for compressible magnetohydrodynamics by Carlier and Campos-Pinto \cite{CarlierCP2024}. These methods, however, struggle with nonsmooth solutions and shocks as they either employ nonconservative variables or lack nonlinear stabilization mechanisms.

In this paper, we present a novel asymptotic-preserving semi-implicit finite element method for weakly compressible flows based on compatible finite element spaces. Our approach employs a semi-implicit time discretization that avoids a Courant-Friedrichs-Lewy (CFL) condition dependent on the sound speed and which results in symmetric positive definite (SPD) linear systems to be solved in each time step. For the spatial discretization, we use compatible finite elements: Raviart-Thomas elements for momentum and discontinuous elements for density, ensuring pointwise conservation of mass. To handle shocks and prevent spurious oscillations, we incorporate an \emph{a posteriori} artificial viscosity limiter based on the MOOD approach by Clain, Diot, and Loub\`{e}re \cite{MOOD, MOODhighorder,MOODorg} and the ideas outlined in \cite{TD17,HybridMHD2}. 
\color{black}
Finally, our method makes use of the celebrated hybridization technique originally introduced by Arnold and Brezzi \cite{ArBr85}, which allows for a more efficient implementation. 
\color{black}

The proposed scheme is asymptotic preserving (AP) in the low Mach number limit, i.e. when the Mach number approaches zero and the density is constant, we obtain a scheme for the incompressible Navier-Stokes equations which has many points in common with the MEEVC scheme by Palha and Gerritsma \cite{MEEVC} and the HDG scheme by Lehrenfeld and Sch\"{o}berl \cite{LeSc16}. However, to the very best of our knowledge, the proposed methodology is completely new in the context of weakly compressible flows and this is also one of the first works that uses \textit{a posteriori} limiting in the context of compatible finite element exterior calculus (FEEC).  

The rest of the paper is organized as follows. In Section~\ref{sec:equations} we recall the equations for weakly compressible flows, discussing a possible viscous regularization and their asymptotic limit when the Mach number approaches zero. In Section~\ref{sec:method} we describe our numerical method. In Section~\ref{sec:results} we validate the proposed scheme with several test cases.

\color{black}
\section{Governing equations}
\color{black}
\label{sec:equations}
We first consider the non-conservative density-momentum-entropy formulation of inviscid weakly compressible isentropic flows
\begin{subequations}
	\begin{align}
		\ddt \density + \Div \mom &= 0, \\
		\ddt \mom + \Div \mathcal{F}(\mom )  &= 0, \\
		\ddt \entropy  + \frac{\mom}{\density}\cdot \Grad \entropy &= 0.
	\end{align}
	\label{eq:compressible}
\end{subequations}
Here $\density$, $\mom$ and $\entropy$ denote density, momentum and specific entropy, respectively. Although the equations can be formulated solely in these three variables, it is useful consider also the velocity $\vel = \mom/ \density$ and the pressure $\pres = \pres(\density, \entropy)$, which is a function of the density and the entropy. With these additional variables, the momentum flux can be defined as 
\begin{equation*}
	\mathcal{F} = \mathcal{F}_{\mom} + \mathcal{F}_{\pres} = \vel \otimes \mom + \pres \, \mathbf{I}, 
\end{equation*}
with $\mathbf{I}$ being the identity tensor. The square of the isentropic speed of sound is defined as 
\begin{equation*}
	c^2 \doteq \frac{\partial \pres}{\partial \density}.
\end{equation*}
Equivalently, it is possible to consider the density as a function of pressure and entropy, i.e. $\density = \density(\pres, \entropy)$. Then the inverse function theorem immediately gives
\begin{equation*}
	\frac{\partial \density}{\partial \pres} = \frac{1}{c^2}.
\end{equation*}
In this work we consider the ideal gas equation of state:
\begin{equation*}
\pres(\density, \entropy) = \density^{\gamma}e^{\entropy/c_v}.
\end{equation*}
Here $c_v$ is the specific heat at constant volume and $\gamma$ is the ratio of specific heats. If not stated otherwise, we choose $c_v = 2.5$ and $\gamma = 1.4$.
\subsection{Viscous regularization}
In this work we consider the following viscous regularization of \eqref{eq:compressible}:
\begin{subequations}
	\begin{align}
		\ddt \density + \Div \mom - \Div(\epsilon_{\density}\Grad \density) &= 0, \\
		\ddt \mom + \Div \mathcal{F}(\mom ) - \Grad \epsilon_{\mom}\Div \mom + \Curl \epsilon_{\mom}\Curl \mom   &= 0, \\
		\ddt \entropy  + \frac{\mom}{\density}\cdot \Grad \entropy - \Div \epsilon_{\entropy}\Grad \entropy &= 0. \label{eq:visc_compressible_S}
	\end{align}
	\label{eq:visc_compressible}
\end{subequations}
Here $\epsilon_Q$ with $Q\in \{ \density, \mom, \entropy\}$ indicates a viscosity, either physical or numerical, which will be specified later. 
\begin{remark}
	When $\epsilon_{\density}$, $\epsilon_{\mom}$ and $\epsilon_{\entropy}$ are constant, this viscous regularization coincides with the \lq\lq monolithic regularization\rq\rq\, discussed by Guermond and Popov in \cite{GuermondViscReg}, but in general they are different. Nevertheless, both regularizations do not have a physical meaning: they serve as a numerical tool to resolve discontinuities and shocks.
\end{remark}
\subsection{Asymptotic limit}
Let us assume now that $\epsilon_{\density} = \epsilon_{\entropy} = 0$ and $\epsilon_{\mom} = \nu$ with $\nu$ being a nonnegative constant.  Then, an asymptotic analysis (see, e.g., \cite{KlaMaj,KlaMaj82,munzMPV,Klein2001}) reveals that when $c^2$ goes to infinity and $\rho$ is constant, the system \eqref{eq:compressible} tends to the incompressible Navier-Stokes equations with the viscosity term in rotational form:
\begin{subequations}
	\begin{align}
		\ddt (\density \vel ) + \Div \mathcal{F}_{\mom}(\density \vel) + \Grad \pres + \mu \Curl \Curl \vel &= 0, \label{eq:NSa}\\
		\Div \vel &= 0, \label{eq:divu}
	\end{align}
	\label{eq:NS}
\end{subequations}
with $\mu = \nu\rho$.
\color{black}
\textcolor{black}{Only for constant viscosity,} the rightmost term in \eqref{eq:NSa} coincides with the usual $ - \mu \Delta \vel$ as a consequence of the following identity:
\begin{equation*}
-\mu \Delta \vel = \mu (\Curl\Curl \vel - \Grad \Div \vel) = \mu  \Curl\Curl \vel.
\end{equation*}
This reformulation of the viscous stress tensor has been put forward for the first time by N\'{e}d\'{e}lec in \cite{Nedelec1982} and has been studied theoretically and numerically for the Stokes and the Navier-Stokes problems, for example, in \cite{DuSaSa2003, DuSaSa2003b, ArFaGo2012, BoFu23, Hanot23, CaCPFa2023}, without claiming completeness. 
\color{black}
\begin{remark}
	We are considering the Laplace formulation of the viscous stress tensor, which is known to yield unphysical solutions in the presence of Navier slip boundary conditions on curved boundaries \cite{ObjectivityNS}. In the case of the rotational formulation of the Laplacian, this problem can be solved by adding an appropriate boundary term proportional to the Weingarten map (see the work of Mitrea and Monniaux \cite{MiMo09}), which is nonzero only on curved boundaries. For simplicity, in this work we consider only flat boundaries when dealing with Navier slip boundary conditions.
\end{remark}
\color{black}

\section{Numerical method} \label{sec:method}

\subsection{Time discretization}
\label{sec:time_disc}
\textcolor{black}{In our semi-implicit scheme the convection of momentum and entropy is treated explicitly, whereas the remaining terms are treated implicitly. For simplicity, we present a low order semi-implicit splitting, see e.g. \cite{CG84,Casulli1990}, with higher order in time achievable via the IMEX methodology, see e.g. \cite{PareschiRusso2000,DLDV18,BDLTV2020,Thomann2020,Thomann2020b,Thomann2022,Thomann2022b,BDLTV2020}. As anticipated in the introduction, we are going to use the density as a function of the pressure and the entropy. Therefore, our time discretization of \eqref{eq:visc_compressible} reads
\begin{subequations}
	\begin{align}
		\density(\pres^{n+1}, \entropy^{n+1})  + \dt \Div \mom^{n+1} &= \density^{n},\label{eq:time-disc_rho} \\
		\mom^{n+1} + \dt \Grad \pres^{n+1} -\dt \Grad(\epsilon_{\mom}\Div\mom^{n+1}) + \dt \Curl (\epsilon_{\mom}\Curl \mom^{n+1}) &= \mom^{n} -\dt\Div\mathcal{F}_{\mom}(\density^n,\mom^n), \label{eq:time-disc_mom}\\
		\entropy^{n+1} - \dt\Div(\epsilon_{\entropy}\Grad \entropy^{n+1}) &= \entropy^{n} - \dt\frac{\mom^{n}}{\density^{n}}\cdot \Grad \entropy^n.\label{eq:time-disc_e}
	\end{align}
	\label{eq:time-disc}
	\end{subequations}
	We solve \eqref{eq:time-disc_rho}-\eqref{eq:time-disc_mom} with the Netwon method. For notational convenience, set $(c^2)^{n+1, l} = c^2(\pres^{n+1, l}, \entropy^{n+1})$, $\mom^* = \mom^n - \dt \Div \mathcal{F}_{\mom}(\density^n,\mom^n)$ and $\density^{n+1, l} = \density(\pres^{n+1, l}, \entropy^{n+1})$, with $\rho^{n+1, 0} = \rho^n$. With these conventions the Newton iteration reads
	\begin{subequations}
	\begin{align}
	\begin{split}
	\frac{1}{(c^2)^{n+1, l}}\pres^{n+1, l+1}  + \dt \Div \mom^{n+1, l+1} &= \density^{n} - \density^{n+1, l} \\
	&+ \frac{1}{(c^2)^{n+1, l}}\pres^{n+1, l},\end{split} \\
	\mom^{n+1, l+1} + \dt \Grad \pres^{n+1, l+1} -\dt \Grad(\epsilon_{\mom}\Div\mom^{n+1, l+1}) + \dt \Curl (\epsilon_{\mom}\Curl \mom^{n+1, l+1})  &= \mom^*.
	\end{align}
	\label{eq:time-disc_Newton}
\end{subequations}
Finally, we update the density:
\begin{equation}
	\density^{n+1} - \Delta t\Div(\epsilon_{\density}\Grad \density^{n+1})= \density^n - \dt \Div \mom^{n+1}.
	\label{eq:density_update}
\end{equation}
\begin{remark}
	Note that the $\Div(\epsilon_{\density}\Grad \density)$ term is neglected in \eqref{eq:time-disc_rho}, but is present in \eqref{eq:density_update}. This is equivalent to a splitting of the density convection and diffusion terms. In this way, the nonlinear system \eqref{eq:time-disc} is considerably simplified.
\end{remark}
}
\subsection{Space discretization}

\subsubsection{Notation}
We first introduce some notation. The $L^2$ scalar product of two scalar-valued or vector-valued functions is indicated by $(\cdot, \cdot)$. Let $\mathcal{T}_h$ be a simplicial triangulation of the domain $\Omega$. A generic element in $\mathcal{T}_h$ is denoted by $T$; its boundary is $\partial T$. A generic facet is denoted by $e$. We denote the $L^2$ product on  $D$ by $\langle \cdot, \cdot\rangle_{D}$ for $D \in \{\partial\Omega, \partial T, e\}$. We define the skeleton of the mesh as $\partial\mathcal{T}_h \doteq \bigcup_{T\in\mathcal{T}}\partial T$, in which each internal facet is counted twice, so that any function on $\partial\mathcal{T}_h$ is \emph{double-valued} on each internal facet. We define the scalar product on $\partial\mathcal{T}_h$ as
\begin{equation*}
	\langle \cdot ,\cdot \rangle_{\partial \mathcal{T}_h} \doteq \sum_{T\in\mathcal{T}_h}\langle \cdot, \cdot \rangle_{\partial T}.
\end{equation*} 
Given a double-valued function $\widehat{v} \in L^2(\partial \mathcal{T}_h)$, we define its average and jump as the elements in $L^2(\partial \mathcal{T}_h)$, such that at a facet $e$ we have 
\begin{align*}
	\ave{\widehat{v}}_{e^{\pm}} &\doteq\begin{cases} \frac{1}{2}(\widehat{v}^+ + \widehat{v}^-) \text{ if $e= \partial T^+\cap \partial T^-$ is an internal facet,}\\ \widehat{v}^+\text{ if $e$ is a boundary facet;}\end{cases}\\
	 \qquad \jump{\widehat{v}}_{e^{\pm}} &\doteq \begin{cases} \textcolor{black}{\pm(\widehat{v}^+ - \widehat{v}^-)}\text{ if $e = \partial T^+ \cap \partial T^-$ is an internal facet,}\\ 0 \text{ if $e$ is a boundary facet.}\end{cases}
\end{align*}
The average is single-valued, i.e. $\ave{\widehat{v}}_{e^+} = \ave{\widehat{v}}_{e^-}$, while the jump is \emph{single-valued up to its sign}, i.e. $\jump{\widehat{v}}_{e^+} = - \jump{\widehat{v}}_{e^-}$.

\textcolor{black}{
 We will consider the following finite element spaces. Let $\Sigma_{r+1}$ be the space of classical continuous Lagrange finite elements of degree $r + 1$ when the dimension is two and the space of N\'{e}d\'{e}lec edge-elements of first kind of degree $r+1$ \cite{Nedelec1} when the dimension is three. We will also indicate by $\RT_{r}$ and $\DG_r$ the spaces of Raviart-Thomas \cite{RT} and discontinuous finite elements respectively. \textcolor{black}{In particular note that $\Curl \Sigma_{r+1} \subset \RT_r$ and $\Div \RT_r = \DG_r$.} For a finite element space $V\in\{\Sigma_{r+1}, \RT_r\}$, we denote by $\mathring{V}$ its subspace with essential boundary conditions. Finally, let $a_{\epsilon}:\DG_r\times \DG_r\to \mathbb{R}$ be the classical symmetric interior-penalty bilinear form associated to the $-\Div(\epsilon\Grad)$ operator introduced by Arnold \cite{SIP}:
\begin{equation*}
a_{\epsilon}(p_h, q_h) \doteq (\epsilon \Grad p_h, \Grad q_h)_{\mathcal{T}_h} + \langle \frac{\zeta \epsilon}{h} \jump{p_h}, \jump{q_h}\rangle_{\partial \mathcal{T}_h}- \langle \ave{\epsilon \Grad p_h}\cdot \nvec, \jump{q_h}\rangle_{\partial\mathcal{T}_h} -  \langle \ave{\epsilon \Grad q_h}\cdot \nvec, \jump{p_h}\rangle_{\partial\mathcal{T}_h}.
\end{equation*}
Here $\zeta$ is a user-defined parameter. In this work we choose $\eta= 40$.}

\subsubsection{Path-conservative DG scheme for the entropy}

\color{black} In eqn. \eqref{eq:visc_compressible_S} the entropy transport equation is used instead of the total energy conservation law, hence we cannot expect shock waves to be correct for large shock Mach numbers, see \cite{toro-book,HouLF94}. This is because non-conservative systems lack a definition of weak solution in the presence of discontinuities and since for the correct computation of shock waves in fluids total energy conservation is mandatory. In this work, we therefore deliberately only consider \textit{weakly compressible} isentropic flows, with Mach numbers ranging from zero to about unity. In general, the appropriate numerical discretization of non-conservative hyperbolic equations still remains a challenge.  
From the theoretical side, a possibile solution to the problem has been proposed by Dal Maso, LeFloch and Murat \cite{DLMtheory}, who introduced a theory (called DLM theory in the following) of weak solutions using paths in phase-space. The DLM theory has inspired Par\'{e}s to develop a a theoretical framework of path-conservative numerical methods \cite{pares2006numerical}. This framework has been used to design Finite Volume methods for non-conservative systems by Par\'{e}s, Castro and collaborators \cite{Castro2006, Castro2007,Munoz2007, Castro2008, CastroPardoPares, Castro2009}. The first path-conservative Discontinuous Galerkin finite element methods have been proposed in \cite{Rhebergen2008,ADERNC} and \cite{USFORCE2}. 

Recall that the velocity $\vel$ is a function of density and momentum, i.e. 
Let $\mathbf{Q} = (\density, \mom)$. Then $\vel = \vel(\mathbf{Q})$ is a function of $\mathbf{Q}$, since $\vel = \mom/ \density$. Let $\Psi(\mathbf{Q}^+, \mathbf{Q}^-)$ be a path in phase space joining $\mathbf{Q}^-$ and $\mathbf{Q}^+$, that is $\Psi(0) = \mathbf{Q}^+$ and $\Psi(1) = \mathbf{Q}^-$. At each mesh interface, we look for a \lq\lq Roe-type normal velocity\rq\rq\, $\widehat{\vel\cdot \nvec}$ satisfying the generalized Rankine-Hugoniot conditions:
\begin{equation}
	\widehat{\vel\cdot \nvec} \jump{S_h} = \int_{0}^1 \vel(\Psi(\mathbf{Q}^+, \mathbf{Q}^-) )\cdot \nvec\frac{\partial \Psi}{\partial s} \dS. \label{eq:gen_RH}
\end{equation}
In this work we choose the segment path $\Psi(s) = (1-s)\mathbf{Q}^+ + s\mathbf{Q}^-$, which yields the following expression for $\widehat{\vel \cdot \nvec}$:
\begin{equation}
\widehat{\vel\cdot \nvec} = \int_0^1 \vel(\Psi(\mathbf{Q}^+, \mathbf{Q}^-))\cdot \nvec \dS. \label{eq:Roe_vel_path}
\end{equation}
The integral in \eqref{eq:Roe_vel_path} can be approximated with a quadrature rule. In our numerical experiments we have noticed that the midpoint rule is sufficient to preserve the accuracy and correctness of the scheme. 
Then, the weak problem associated to each time-step reads as follows.
\begin{weakproblem}
	Find $\entropyh^{n+1}$ such that 
	\begin{equation*}
	(\entropyh^{n+1}, \entropytesth) + \dt a_{\epsilon_\entropy}(\entropy^{n+1}_h, \entropytest_h) = (\entropyh^n, \entropytesth) - \dt \left( \frac{\momh^n}{\densityh^n}\cdot \Grad \entropyh^{n}, \entropytesth \right) + \frac{\dt}{2} \langle (\widehat{\vel\cdot \nvec} - s_{\max}) \jump{\entropyh^n}, \entropytesth\rangle_{\partial \mathcal{T}_h}
	\end{equation*}	
	for each $\entropytesth\in \DG_r$. 
\end{weakproblem} 
Here $s_{\max} \doteq \max\left( 2\left\lvert\frac{\momh^+}{\densityh^+}\cdot \nvec \right\rvert , 2\left\lvert\frac{\momh^-}{\densityh^-}\cdot \nvec\right\rvert\right)$. Note that the matrix associated to this linear system is symmetric positive definite.

\subsubsection{Convection of the momentum}
{For the convection of the momentum we employ a standard DG discretization. Let $\mom^*_h$ be the solution of 
\begin{equation}
	(\momh^*, \veltesth) = (\momh^n, \veltesth) - \dt (\mathcal{F}_{\mom}(\rho^n,\momh^n), \Grad_h \veltesth) + \langle \widehat{\mathcal{F}}_{\mom}(\densityh^n,\momh^n)\nvec , \veltesth\rangle_{\partial\mathcal{T}_h}. \label{eq:DG_mom}
\end{equation}
At an interface $e = \partial T^+\cap \partial T^-$, in order to discretize the nonlinear convective terms we use a Ducros-type numerical flux with a dissipative term:
\begin{equation}
	\widehat{\mathcal{F}}(\densityh,\momh)\nvec = \momh\cdot\nvec  \ave{{\momh}/{\densityh}} + \frac{1}{2}s_{\max} \, \jump{ \momh},
\end{equation}
with $s_{\max}$ is defined as in the previous section. 

\subsubsection{Mixed finite element discretization of the momentum-pressure system}
We introduce the \lq\lq momentum vorticity\rq\rq\, $\vor = \epsilon_{\mom}\Curl \mom$. We approximate $\vor$ with $\Sigma_{r+1}$,  $\mom$ with Raviart-Thomas elements $\RT_{r}$ \cite{RT} and $\pres$ with element-wise discontinuous polynomials $\DG_r$. Then, at each Newton iteration we solve the following linear system.
\begin{weakproblem}
	Find $\vorh\in \Sigma_{r+1}$, $\momh^{n+1, l+1}\in \RT_{r}$, $\presh^{n+1, l+1}\in \DG_r$ satisfying 
	\begin{subequations}
		\begin{align}
			\left(\frac{1}{(c^2)^{n+1, l}}\pres^{n+1, l+1}_h, \prestesth \right)  + \dt (\Div \momh^{n+1, l+1}, \prestesth) &= \left(\density^{n}_h - \density^{n+1, l} + \frac{1}{(c^2)^{n+1,l}}\presh^{n+1, l}, \prestesth\right), \label{eq:discrete_Newton_rho}\\
		\begin{split}	(\momh^{n+1, l+1}, \veltesth) - \dt( \presh^{n+1, l+1}, \Div\veltesth) & \\
		+ \dt(\epsilon_{\mom}\Div \mom_h^{n+1, l+1}, \Div \veltesth) + \dt(\Curl \vorh^{n+1, l+1}, \veltesth) &= (\mom^*_h, \veltesth) - \langle \overline{p}, \veltesth\cdot \nvec \rangle_{\partial \Omega}, \end{split}\label{eq:discrete_Newton_m}\\
			\left(\frac{1}{\epsilon_{\mom}}\vorh^{n+1, l+1}, \vortesth\right) - (\mom_h^{n+1, l+1}, \Curl \vortesth)& = \langle \nvec \times \overline{\mom}, \vortesth\rangle_{\partial \Omega}, \label{eq:discrete_Newton_vor}
		\end{align}
		\label{eq:discrete_Newton}
	\end{subequations}
	for each $\prestesth\in \DG_r$, $\veltesth\in \RT_{r}$ and $\vortesth\in \Sigma_{r+1}$. 
	\end{weakproblem}
\begin{remark}
	\label{rmk:bc}
	This choice of spaces corresponds to \lq\lq outflow\rq\rq\, boundary conditions:
	\begin{equation*}
	(\pres - \epsilon_{\mom}\Div \mom)\rvert_{\partial \Omega} = \overline{p}, \qquad \nvec \times \mom\rvert_{ \partial \Omega} = \nvec \times \overline{\mom}.
	\end{equation*}
	If we choose the spaces $\mathring{\Sigma}_{r+1}$, $\mathring{\RT}_r$ in place of $\Sigma_{r+1}$ and $\RT_r$, we obtain the nonstandard \lq\lq slip\rq\rq\, boundary conditions:
	\begin{equation*}
	\mom\rvert_{ \partial \Omega} \cdot \nvec = 0, \qquad \nvec \times(\Curl \mom)\rvert_{ \partial \Omega} = 0.
	\end{equation*}
	See Mitrea and Monniaux \cite{MiMo09} for a discussion on how these boundary conditions relate to the standard Navier slip ones. Finally, the choice $\Sigma_{r+1}$ and $\mathring{RT}_r$ yields Dirichlet boundary  conditions (which include both those commonly referred as \lq\lq wall\rq\rq\, and \lq\lq inflow\rq\rq):
	\begin{equation*}
	\mom\rvert_{ \partial \Omega} = \overline{\mom}. 
	\end{equation*}
	In this latter case, the inclusion $\Curl \Sigma_{r+1}\subset \mathring{\RT}_r$ is false. This issue is the source of the subpotimal convergence rate of $\omega_h$ and $p_h$ shown by Arnold, Falk and Gopalakrishnan in \cite{ArFaGo2012}.
\end{remark}
When the Newton method has converged, we update the density solving the problem 
\begin{equation}
( \density^{n+1}_h, q_h) + \Delta t a_{\epsilon_{\density}}(\density^{n+1}_h, q_h) = ( \density^n_h, q_h) - \dt ( \Div \momh^{n+1}, q_h).
\label{eq:update_rho}
\end{equation}
We can immediately deduce the following conservation properties of our scheme.
\begin{theorem}
Assume that $\momh\cdot \nvec$ and $\nvec \times \vorh$ vanish on the boundary of $\Omega$ and $\epsilon_{\density} = 0$. If the Newton iteration \eqref{eq:discrete_Newton} converges, then the resulting scheme conserves mass locally and momentum globally, that is:
\begin{align}
\int_T \density^{n+1}_h \dx &= \int_T \density^{n}_h \dx + \dt \int_{\partial T}\mom^{n+1}_h \cdot \nvec \dS, \qquad \forall T\in \mathcal{T}_h, \\
\int_{\Omega}\mom^{n+1}_h \dx &= \int_{\Omega}\mom^{n}_h\dx.
\end{align}
Furthermore, if either $\epsilon_{\density} = 0$ or $\Grad \rho \cdot \nvec= 0$ on $\partial\Omega$,  mass is conserved also globally:
\begin{equation*}
\int_{\Omega}\rho^{n+1}\dx = \int_{\Omega}\rho^n \dx.
\end{equation*}
\end{theorem}
\begin{proof}
When $\epsilon_{\density} = 0$, equation \eqref{eq:update_rho} reduces to the simple update $\rho^{n+1}_h = \rho^n_h - \Delta t \Div \momh^{n+1}$, since $\Div \RT_{r} = \DG_r$. The first property follows then from the Gauss theorem. To prove the second property, take $\veltesth = \mathbf{e}_i$ in \eqref{eq:DG_mom} and \eqref{eq:discrete_Newton_m}. Clearly, the terms involving derivatives of $\mathbf{e}_i$ vanish. It remains to show that $\langle\widehat{\mathcal{F}}_{\mom}(\density^n_h, \mom^n_h)\cdot \nvec , \mathbf{e}_i\rangle_{\mathcal{T}_h}$ and $(\Curl \vor_h^{n+1}, \mathbf{e}_i)$ vanish as well. The first term is zero since $\mathbf{e}_i$ is continuous and $\widehat{\mathcal{F}}_{\mom}(\rho^n_h, \mom^n_h)\cdot \nvec$ is single-valued at internal facets and vanishes on the boundary. For the second one, the claim follows from the simple computation 
\begin{equation*}
\begin{split}
 (\Curl \vorh^{n+1}, \mathbf{e}_i) &= \int_{\Omega}\Curl \vorh^{n+1}\cdot \mathbf{e}_i \dx\\ &= \int_{\Omega} \Div (\vorh^{n+1} \times \mathbf{e}_i) \dx\\  &= \int_{\partial \Omega}(\vorh^{n+1}\times \mathbf{e}_i)\cdot \nvec \dS\\
 &= \int_{\partial \Omega}(\nvec \times\vorh^{n+1} )\cdot\mathbf{e}_i \dS\\
 &= 0.
 \end{split}
\end{equation*} 
In the last step we have used the fact that $\nvec \times \vorh$ vanishes on the boundary. Finally, if $\epsilon_{\density} = 0$, the global conservation of mass follows from the local one since $\momh\cdot \nvec$ is single-valued at internal facets and vanishes on the boundary. On the other hand, if $\epsilon_{\density} \neq 0$ and $\Grad \density\cdot \nvec = 0$ on $\partial \Omega$, the claim follows by taking $q_h = 1$ in \eqref{eq:update_rho} since $a_{\epsilon_{\density}}(\rho^{n+1}_h, 1) = 0$.
\end{proof}

\subsubsection{Efficient decoupling of the vorticity from the pressure}
Taking $\veltesth = \Curl \vortest_h$ with $\vortest_h\in \Sigma_{r+1}$ in \eqref{eq:discrete_Newton_m} and using $\Div\Curl = 0$, we obtain 
\begin{equation}
(\momh^{n+1, l+1}, \Curl \vortest_h) + \dt(\Curl \vorh^{n+1, l+1}, \Curl \vortest_h) = (\mom^*_h, \Curl \vortest_h) - \langle\overline{p}, \Curl \vortest_h \cdot \nvec\rangle_{\partial \Omega}. \label{eq:reduced_mom}
\end{equation}
Taking the sum of \eqref{eq:reduced_mom} with \eqref{eq:discrete_Newton_vor}, we obtain a single decoupled equation for the vorticity. 
\begin{weakproblem}
	Find $\vorh^{n+1, l+1}\in \Sigma_{r+1}$ satisfying 
\begin{equation}
\left(\frac{1}{\epsilon_{\mom}}\vorh^{n+1, l+1}, \vortesth\right) + \dt(\Curl \vorh^{n+1, l+1}, \Curl \vortest_h) = (\mom^*_h, \Curl \vortest_h) - \langle\overline{p}, \Curl \vortest_h \cdot \nvec\rangle_{\partial \Omega} + \langle\nvec \times \overline{\mom}, \vortesth\rangle_{\partial \Omega},
\label{eq:vorticity_dec}
\end{equation}
for each $\vortest_h\in \Sigma_{r+1}$. 
\end{weakproblem}
Note that the matrix associated to this linear system is symmetric positive definite. Once we have computed $\vorh^{n+1, l+1}$, we can compute momentum and pressure as follows. 
\begin{weakproblem}
	Find $\momh^{n+1, l+1}\in\RT_r$ and $\presh^{n+1, l+1}\in \DG_r$ satisfying
	\begin{subequations}
		\begin{align}
		\left(\frac{1}{(c^2)^{n+1, l}}\pres^{n+1, l+1}_h, \prestesth \right)  + \dt (\Div \momh^{n+1, l+1}, \prestesth) &= \left(\density^{n}_h - \density^{n+1, l} + \frac{1}{(c^2)^{n+1,l}}\presh^{n+1, l}, \prestesth\right), \\
		\begin{split}(\momh^{n+1, l+1}, \veltesth) - \dt( \presh^{n+1, l+1}, \Div\veltesth) &\\ + \dt(\epsilon_{\mom}\Div \mom_h^{n+1, l+1}, \Div \veltesth) &= (\mom^*_h, \veltesth)- \langle \overline{p}, \veltesth\cdot \nvec \rangle_{\partial \Omega} - \dt(\Curl \vorh^{n+1, l+1}, \veltesth),\end{split}
		\end{align}
		\label{eq:discrete_Newton_pm}
	\end{subequations}
for each $\veltesth\in \RT_r$ and $\prestesth\in \DG_r$.
\end{weakproblem} 
\color{black}
\begin{remark}
The same splitting procedure can be applied in the case of nonstandard slip boundary conditions since $\Curl \mathring{\Sigma}_{r+1} \subset \mathring{\RT}_r$. On the other side, the case of Dirichlet boundary condition is more delicate, since we cannot take $\veltesth = \Curl \vortest_h$ due to $\Curl \Sigma_{r+1} \not\subset \mathring{RT}_r$, as already explained in Remark \ref{rmk:bc}. We can circumvent this obstacle with the following trick. Instead of imposing the boundary condition $\mom\rvert_{ \partial \Omega} \cdot \nvec = 0$ essentially, we impose it via a Lagrange multiplier $\widehat{p}_h$ belonging to the space of discontinuous polynomials on the boundary:
\begin{equation*}
\widehat{M}_h = \{ \widehat{q}_h \in L^2(\partial \Omega) \mid \widehat{q}_h\rvert_{ e}\in \mathcal{P}_r(e)\, \forall e\subset \partial \Omega\}.
\end{equation*}
Using this space, we can rewrite equation \eqref{eq:discrete_Newton_m} as the following equivalent system: find $(\momh^{n+1, l+1}, \widehat{p}_h)\in \RT_r \times \widehat{M}_h$ satisfying
\begin{subequations}
	\begin{align}
	\begin{split}	(\momh^{n+1, l+1}, \veltesth) - \dt( \presh^{n+1, l+1}, \Div\veltesth) & \\
	+ \dt(\epsilon_{\mom}\Div \mom_h^{n+1, l+1}, \Div \veltesth) + \dt(\Curl \vorh^{n+1, l+1}, \veltesth) + \langle \widehat{p}_h, \veltesth\cdot \nvec \rangle_{\partial \Omega} &= (\mom^*_h, \veltesth) , \end{split} \label{eq:discrete_Newton_m_LM}\\
	\langle \mom_h^{n+1, l+1}\cdot \nvec, \widehat{q}_h\rangle_{\partial \Omega} &= 0,
	\end{align}
\end{subequations}
for each $(\veltesth, \widehat{q}_h)\in \RT_r \times \widehat{M}_h$. Now it is possible to take $\veltesth = \Curl \vortest_h$ in \eqref{eq:discrete_Newton_m_LM}, obtaining the following equation for $\vor^{n+1,l+1}_h$:
\begin{equation}
\left(\frac{1}{\epsilon_{\mom}}\vorh^{n+1, l+1}, \vortesth\right) + \dt(\Curl \vorh^{n+1, l+1}, \Curl \vortest_h) + \langle\widehat{p}_h, \Curl \vortest_h \cdot \nvec\rangle_{\partial \Omega} = (\mom^*_h, \Curl \vortest_h)  + \langle\nvec \times \overline{\mom}, \vortesth\rangle_{\partial \Omega},
\label{eq:vorticity_dec_LM}
\end{equation}
\color{black}
The presence of $\widehat{p}$ makes equation \eqref{eq:vorticity_dec_LM} an underdetermined problem. To recover well-posedness, we take either $\widehat{p}_h = \overline{p}$, when the boundary pressure is known, or $\widehat{p}_h = p^n_h\rvert_{ \partial \Omega}$. Then, the \textcolor{black}{momentum} and the pressure are obtained as in the other cases. 
\end{remark}
\color{black}
\paragraph{Linear algebra and hybridization}
We briefly comment on the algebraic structure of system \eqref{eq:discrete_Newton_pm} and the possible solution strategy. Let $M_{c^2}^{\pres}$ and $M^{\mom}$ the (weighted) mass matrices associated to $\pres$ and $\mom$ respectively, and let $D$ and $S$ be the div and div-div matrices. Then, the matrix $G$ associated to \eqref{eq:discrete_Newton_pm} reads
\begin{equation*}
G = \begin{pmatrix} M_{c^2}^{\pres} & D\\ -D^T & M^{\mom} + S\end{pmatrix}.
\end{equation*}
Now, since $M_{c^2}^{\pres}$ is block diagonal, it can be inverted cheaply and we can consider its Schur complement $\widetilde{G}_{\mom} = M + S + D^T(M_{c^2}^{\pres})^{-1}D$ which is symmetric and positive definite. However, when $c^2 \to \infty$,  the matrix $M_{c^2}^{\pres}$ tends to $0$ and $\widetilde{G}_{\mom}$ becomes ill-conditioned, so this strategy is impractical for large values of $c^2$. 
\color{black}
As a remedy, we introduce hybridization, that is, we break the normal continuity of the space $\RT_r$ and we enforce it via a Lagrange multiplier. \color{black} Let $\widetilde{\RT}_r$ be the \lq\lq broken version\rq\rq\, of $\RT_r$, that is $\widetilde{\RT}_r$ is the space of $L^2$ functions on $\Omega$ such that the restriction on each element is a Raviart-Thomas polynomial of degree $r$. Now, given a facet $e$, ${\veltesth}\rvert_{ e}\cdot \nvec\in \mathcal{P}_{r-1}(e)$ for $\veltesth\in \RT_r$, where $\mathcal{P}_{r-1}(e)$ is the space of polynomials of degree $r-1$ on $e$ (for a proof, see Proposition 2.3.3 in \cite{BoBrFoBook}). It follows that the normal continuity of $\mom_h$ can be imposed via a Lagrange multiplier ${\lambda}_h$ in the space 
\begin{equation*}
M_h \doteq \prod_{e\in{\mathcal{E}}_h}\mathcal{P}_{r-1}(e).
\end{equation*}
Define also the space $\mathring{M}_{\overline{\pres},h}$ as the subspace of $M_h$ \lq\lq with boundary conditions\rq\rq:
\begin{equation*}
\mathring{M}_{\overline{\pres},h} \doteq \{ {\lambda}_h \in M_h\mid \langle {\lambda}_h, {\xi}_h\rangle_e = \langle \overline{p}, {\xi}_h\rangle_e\,\forall {\xi}\in\mathcal{P}_r(e)\,\forall e\subset\partial\Omega\}.
\end{equation*}
\color{black}
The hybridized formulation with outflow boundary conditions reads:
\begin{weakproblem}
	Find $\momh^{n+1, l+1}\in\widetilde{\RT}_r$, $\presh^{n+1, l+1}\in \DG_r$ and $\lambda_h\in \mathring{M}_{\overline{\pres},h}$ satisfying
	\begin{subequations}
		\begin{align}
		\left(\frac{1}{(c^2)^{n+1, l}}\pres^{n+1, l+1}_h, \prestesth \right)  + \dt (\Div \momh^{n+1, l+1}, \prestesth) &= \left(\density^{n}_h - \density^{n+1, l} + \frac{1}{(c^2)^{n+1,l}}\presh^{n+1, l}, \prestesth\right), \\
		\begin{split}(\momh^{n+1, l+1}, \veltesth) - \dt( \presh^{n+1, l+1}, \Div\veltesth) \\+ \dt(\epsilon_{\mom}\Div \mom_h^{n+1, l+1}, \Div \veltesth) + \langle \lambda_h, \veltesth \cdot \nvec\rangle_{\partial \mathcal{T}_h} &= (\mom^*_h, \veltesth) - \dt(\Curl \vorh^{n+1, l+1}, \veltesth),\end{split}\\ 
		\langle \momh^{n+1, l+1}\cdot \nvec, \xi_h\rangle_{\partial \mathcal{T}_h} &= 0,
		\end{align}
		\label{eq:discrete_Newton_pm_h}
	\end{subequations}
	for each $\veltesth\in \widetilde{\RT}_r$, $\prestesth\in \DG_r$ and $\xi_h \in \mathring{M}_{0,h}$.
\end{weakproblem} 
\begin{remark}
Note that in the hybridized formulation \eqref{eq:discrete_Newton_pm_h}, the boundary condition $(\pres - \epsilon_{\mom}\Div \mom)\rvert_{\partial \Omega} = \overline{p}$ is imposed essentially, rather than naturally. Then, as it is customary in this situation, we consider the space $\mathring{M}_{0,h}$ with homogeneous boundary conditions for the test functions.
\end{remark}
\begin{remark}
	In the case of either Dirichlet or slip boundary conditions, the space $\mathring{M}_{\overline{p},h}$ must be replaced with $M_h$.
\end{remark}
The resulting matrix has the block structure 
\begin{equation*}
\widetilde{G}= \begin{pmatrix} M^{\pres}_{c^2} & \widetilde{D} & 0 \\
\widetilde{D}^T & -(\widetilde{M}^{\mom} + S) & \widetilde{B}^T \\
0 & \widetilde{B}& 0 \end{pmatrix}.
\end{equation*}
Defining the matrices
\begin{equation*}
 L \doteq  \begin{pmatrix}  M^{\pres}_{c^2} & \widetilde{D} \\ \widetilde{D}^T & -(\widetilde{M}^{\mom}+S) \end{pmatrix}, \quad C \doteq \begin{pmatrix} 0 & I\end{pmatrix},
\end{equation*}
with $I$ being the identity matrix, we can rewrite $\widetilde{G}$ as 
\begin{equation*}
\widetilde{G}= \begin{pmatrix} L & C^T\widetilde{B}^T \\ 
\widetilde{B}C & 0\end{pmatrix}.
\end{equation*}
In the parlance of Cockburn, Gopalakrishnan and Lazarov \cite{CoGoLa09}, $L$ is the \lq\lq local solver\rq\rq\, matrix and it is block-diagonal. Moreover, $L$ remains invertible even when $c^2\to \infty$. We can then consider its Schur complement $\widetilde{B}C L^{-1}C^T\widetilde{B}^T$, which is symmetric and negative definite, as we now show.
\begin{lemma}
	The matrix $\widetilde{B}C E^{-1}C^T\widetilde{B}^T$ is symmetric negative definite.
\end{lemma} 
\begin{proof}
	In this proof, given a generic finite element function $v$ in a finite element space $V$, we denote by $\vec{v}$ the associated vector of degrees of freedom.	The symmetry is evident since $L$ is symmetric. We prove the negative definiteness. Given $\vec{\xi}$, let $(\vec{q}, \vec{\mathbf{v}}) = -L^{-1}C\widetilde{B}^T \vec{\lambda}^T$, that is
	\begin{equation*}
	\begin{split}
	-(\widetilde{M}^{\mom}+S)\vec{\mathbf{v}} + \widetilde{D}^T\vec{q} &= -\widetilde{B}^T \vec{\xi}, \\
	M^p_{c^2}\vec{q} + \widetilde{D}\vec{\mathbf{v}} &= 0.
	\end{split}
	\end{equation*}
	In particular, we remark that $(\vec{q}, \vec{\mathbf{v}}) = 0$ if and only if $\vec{\xi}= 0$, due to the invertibility of $L$. 
	Then, we have
	\begin{align*}
	\vec{\xi}^T \cdot  \widetilde{B}C L^{-1}C^T\widetilde{B}^T\vec{\xi} &= \vec{\xi}^T \cdot  \widetilde{B}C \begin{pmatrix} \vec{q} \\ \vec{\mathbf{v}} \end{pmatrix}\\
	&= - \vec{\mathbf{v}}\cdot  (\widetilde{M}^{\mom} + S) \vec{\mathbf{v}} + \vec{q}\cdot \widetilde{D} \vec{\mathbf{v}} \\
	&= - \vec{\mom} \cdot (\widetilde{M}^{\mom} + S) \vec{\mathbf{v}}  - \vec{q}\cdot M^{\pres}_{c^2} \vec{q}.
	\end{align*}
	The claim follows from the positive-definiteness of $\widetilde{M}^{\mom} + S$ and and $M^{p}_{c^2}$.
\end{proof}
\color{black}
\begin{remark}
We remark that hybridization of \eqref{eq:vorticity_dec} is possible (see \cite{HybridFEEC}), but it is difficult to implement and does not yield a significant advantage, since the matrix associated to \eqref{eq:vorticity_dec} is already symmetric positive definite.
\end{remark}
\color{black}

\subsection{A posteriori limiting via artificial viscosity}
\color{black}We now discuss the choices of $\epsilon_{Q}$ with $Q\in \{\density, \mom, \entropy\}$. We follow the \emph{a posteriori} MOOD concept originally introduced by Clain, Diot and Loub\`{e}re \cite{MOOD, MOODhighorder, MOODorg} for finite volume methods. This idea has been successfully applied to high order DG methods \cite{DGlimiter1, DGlimiter2, DGlimiter3} to construct \textit{a posteriori} subcell limiters, and for staggered semi-implicit DG methods it has been successfully used by Tavelli and Dumbser in \cite{TD17}. Recently, we have shown that the MOOD paradigm is effective also for compatible finite elements \cite{HybridMHD2}. The MOOD strategy consists of three steps:
	\begin{enumerate}
		\item Computation of a so-called \textit{candidate solution} at time $n+1$ without the use of any limiting and/or artificial viscosity;
		\item Detection of troubled cells by violation of numerical and/or physical admissibility criteria;
		\item Re-computation of the solution with limiting/artificial viscosity on the troubled cells.
	\end{enumerate}
	Following \cite{TD17}, our detection criterion is based on a discrete maximum principle. Let $W$ be a scalar function depeding possibly on $\{\density, \mom, \entropy\}$. Then, we say that $W$ satisfies the relaxed discrete maximum principle on the element $T\in\mathcal{T}_h$ if 
\begin{equation}
\min_{\mathbf{y} \in \mathcal{N}(T)}W(\mathbf{y}, t^n) - \delta_T \leq W(\mathbf{x}, t^n) \leq \max_{\mathbf{y}\in\mathcal{N}(T)}W(\mathbf{y}, t^n) + \delta_T, \qquad \forall \mathbf{x} \in T.
\end{equation}
Here $\mathcal{N}(T)$ is the set made of the Voronoi neighbors of $T$ and $T$ itself, and $\delta_T$ is a relaxation parameter, which is defined as follows 
\begin{equation}
\delta_T = \max\left(\delta_0, \eta\left(\max_{\mathbf{y}\in\mathcal{N}(T)}W(\mathbf{y}, t^n) -  \min_{\mathbf{y} \in \mathcal{N}(T)}W(\mathbf{y}, t^n)\right)\right)
\label{eq:DMP}
\end{equation}
with $\delta_0$ and $\eta$ user-defined parameters. In this work we choose $\delta_0 = 10^{-4}$ and $\eta =10^{-3}$. Then we set 
\begin{equation}
{\epsilon_Q}\rvert_{ T} = \begin{cases} \frac{1}{2}h s_{\max} \text{ if \eqref{eq:DMP} is violated on $T$,} \\ \overline{\epsilon} \text{ otherwise,}	\end{cases}
\label{eq:def_av}
\end{equation}
with $s_{\max}$ an estimate of the maximum wavespeed in the full system and $\overline{\epsilon}$ being a small value to avoid division by zero. 
\subsection{Summary of the algorithm}
The final algorithm can be summarized in the following steps:
\begin{enumerate}
	\item Compute $\entropy^{n+1}_h$ via the path-conservative DG scheme and $\epsilon_S = 0$;
	\item Compute $\epsilon_S$ as in \eqref{eq:def_av} with $W = \entropy$;
	\item Repeat Step 1 with the new $\epsilon_S$;
	\item Compute $\momh^{n+1}$, $p^{n+1}_h$ and $\density_h^{n+1}$ with the Newton iteration \eqref{eq:discrete_Newton} with $\epsilon_{\density} = \epsilon_{\mom} = 0$. In particular at each iteration, we solve \eqref{eq:vorticity_dec} and \eqref{eq:discrete_Newton_pm_h}.
	\item Compute $\epsilon_{\rho}$ as in \eqref{eq:def_av} with $W= \density$;
	\item Repeat Step 4 with $\epsilon_{\mom} = \epsilon_{\rho}$ with $\epsilon_{\rho}$ computed in the previous Step.
\end{enumerate}
\subsection{Incompressible limit scheme}
In this section, we describe a scheme for the incompressible Navier-Stokes equations with constant density and viscosity. Moreover, we show that this scheme is the limit of the method introduced in the previous sections when the density is constant  and the Mach number goes to zero.}
 In particular, let $\velh^* = \momh^*/\rho$ where $\momh^*$ is defined in \eqref{eq:DG_mom}. Then, the vorticity, the velocity and the pressure at time $t^{n+1}$ are obtained solving the following problem.
\begin{weakproblem}
	Find $(\vorh^{n+1},\velh^{n+1}, \presh^{n+1})\in (\Sigma_{r+1}\times \RT_r\times \DG_r)$ satisfying 
	\begin{subequations}
		\label{eq:fully_discrete}
		\begin{align}
			\left(\frac{1}{\mu}\vorh^{n+1}, \vortesth\right) - (\velh^{n+1}, \Curl \vortesth) & = \langle \overline{\velh}\times\nvec, \vortesth\rangle_{\partial \Omega}, \label{eq:vor_smi}\\ 
			(\rho \velh^{n+1}, \veltesth) + \dt (\Curl \vorh^{n+1}, \veltesth) - \dt (\presh^{n+1}, \Div \veltesth) & = (\rho \velh^*, \veltesth) - \dt\langle \overline{\pres}, \veltesth\cdot\nvec\rangle_{\partial\Omega}, \label{eq:vel}\\
			(\Div \velh^{n+1}, \prestesth) & = 0 \label{eq:pres}, 
		\end{align}
	\end{subequations}
for each $(\vortesth, \veltesth, \prestesth) \in  (\Sigma_{r+1}\times \RT_r\times \DG_r)$. 
\end{weakproblem}
	We prove now that in the zero Mach number limit, the Newton method \eqref{eq:discrete_Newton} converges in one iteration and its solution coincides with the solution of \eqref{eq:fully_discrete}.
\begin{theorem}
	Assume that $\rho_h^{n} = \rho$ is constant in space and $\frac{1}{(c^2)^{n+1, 0}} = 0$. Assume moreover that $\epsilon_{\mom} = \mu/\rho$ with $\mu$ being a constant. Then, let $\mom^{n+1, 1}_h$, $\vorh^{n+1, 1}$ and $\presh^{n+1, 1}$ be solutions of \eqref{eq:discrete_Newton} with $l =0$. Then $\vel_h^{n+1} = \frac{\momh^{n+1, 1}}{\rho}$, $\vorh^n = \vorh^{n+1, 1}$ and $\pres^{n+1} = \pres^{n+1, 1}$ solve \eqref{eq:fully_discrete}. 
\end{theorem}
\begin{proof}
	First, note that $\vel^{n+1}_h$ belongs to $\RT_r$ since $\momh^{n+1, 1}\in\RT_r$ and $\rho$ is a constant. Under the assumptions of the Theorem and $\rho^{n+1, 0} = \rho^n_h$ we obtain that \eqref{eq:discrete_Newton_rho} reduces to 
\begin{equation*}
\Delta t (\Div (\rho\vel^{n+1}_h), q_h) = 0,
\end{equation*}
which is equivalent to \eqref{eq:pres} since $\rho$ is constant. Moreover, this implies that the div-div term in \eqref{eq:discrete_Newton_m} vanishes, and therefore \eqref{eq:discrete_Newton} reduces to \eqref{eq:vel}.  Similarly, taking $\epsilon_{\mom} = \mu/\rho$, \eqref{eq:discrete_Newton_vor} reduces to 
\begin{equation*}
\left( \frac{\rho}{\mu}\vor^{n+1}_h, \vortesth\right) - (\rho\vel^{n+1}_h, \Curl \vortesth) = \langle \rho\overline{\vel}\times \nvec, \vortesth\rangle_{\partial \Omega},
\end{equation*}
which is \eqref{eq:vor_smi} multiplied by $\rho$, concluding the proof.
\end{proof}
\color{black}
The scheme \eqref{eq:fully_discrete} is novel, as the considered combination of the reformulation of the equations, the compatible finite element spaces, the semi-implicit time discretization and the employed hybridization techniques has never appeared in the literature before. Nevertheless, our method has many analogies with some existing schemes for the incompressible Navier-Stokes equations, which we now highlight. We confine ourselves to $H(\mathrm{ div})$-based methods for unsteady incompressible flows that preserve $\Div \vel = 0$ exactly pointwise everywhere. In particular, we neglect discontinuous Galerkin methods that achieve this property via postprocessing, such as those devised by Cockburn and collaborators \cite{CoKaSc2004, CoKaSc2006, NgPeCo2011}. We discuss the following ingredients: 
\begin{itemize}
	\item Spatial discretization of the convective term;
	\item Spatial discretization of the viscous term; 
	\item Time discretization.
\end{itemize}
\paragraph{Convective term}
Our DG-based discretization of the convective term coincides with the one proposed by Guzm\'{a}n et al. in \cite{GuShuSe2017}. An alternative would be rewriting the convective term using Lamb's identity:
\begin{equation*}
\Div(\rho\vel\otimes\vel) = \rho(\vel\cdot\Grad)\vel = \rho(\Curl\vel)\times\vel + \frac{1}{2}\rho\Grad(\vel\cdot\vel)
\end{equation*}
Then, the quadratic term $\frac{1}{2}(\vel\cdot\vel)$ is incorporated in the pressure variable. This strategy has been used in \cite{MEEVC,ZPGR22,MEEVC24,Gawlik20, Hanot23}. Another alternative discretization of the convective term recently appeared in \cite{CaCPFa2023}. The weak formulation of the convective term is rewritten as follows assuming vanishing boundary conditions for the velocity:
\begin{equation*}
\int_{\Omega}\rho(\vel\cdot\Grad)\vel\cdot\veltesth\dx = \frac{1}{2}\int_{\Omega}\rho(\vel\cdot\Grad)\vel\cdot\veltesth\dx - \frac{1}{2}\int_{\Omega}\rho(\vel\cdot\Grad)\veltesth\cdot\vel\dx
\end{equation*}
The term on the right-hand side is then discretized using appropriate projection operators. The resulting methods are conforming (no jump terms appear), however no upwinding is present. Therefore implicit timestepping is mandatory, leading to the solution of a nonlinear nonsymmetric system at each time step. On the other side, DG methods allow dissipative upwinding and therefore are stable also with explicit time discretizations. Finally, we mention the nondissipative upwind DG method of Natale and Cotter \cite{NaCo2018}. The convection term is treated as a discrete Lie derivative as defined by Heumann et al. in \cite{HeuHipPag2016}:
\begin{equation}
(\rho\velh, \Curl(\velh\times \veltesth) - \velh\Div \veltesth)_{\mathcal{T}_h} + \langle \nvec \times \rho\velh^{\mathrm{upw}}, \jump{\vel\times\veltesth}\rangle_{\partial \mathcal{T}_h}.
\label{eq:Lie}
\end{equation}
Note that when $\veltesth = \velh$, \eqref{eq:Lie} vanishes. As a consequence, this spatial discretization, when coupled with a midpoint rule in time, conserves the energy, which is remarkable for an upwind scheme. \textcolor{black}{Finally, we remark that many of the alternatives mentioned here assume that the density is constant, which is not our case in general. See Gawlik and Gay-Balmaz \cite{Gawlik20} for a rotational form of the advection term in the case of nonconstant density.}

\paragraph{Viscous term}
As an alternative to the vorticity-based reformulation of the viscous term used also in \cite{MEEVC, ZPGR22,Hanot23, CaCPFa2023, MEEVC24}, it is possible to use hybridizable discontinuous Galerkin methods (HDG), as done, for example, in \cite{LeSc16} or discontinuous Galerkin (DG) \cite{Fu19}. 

\paragraph{Time discretization}
Our simple semi-implicit time discretization necessitates of the solution of only symmetric positive definite linear systems. Moreover, the implicit treatment of the viscous term avoids a quadratic CFL restriction of the time step size. On the other sidem, all the methods \cite{GuShuSe2017, MEEVC, NaCo2018, Gawlik20, ZPGR22,Hanot23, CaCPFa2023, MEEVC24} treat the convection term implicitly, needing the solution of a \textit{nonsymmetric} system at each time step. The only scheme that avoids the solution of a nonsymmetric system is the semi-implicit scheme proposed by Lehrenfeld and Sch\"{o}berl \cite{LeSc16}, which circumvents the problem via a pseudo-implicit approach. We mention also the scheme of Fu \cite{Fu19}, which treats the both the convective term and viscous term explicitly. This latter method is very cheap, since it requires only the solution of one single symmetric linear system at each time iteration, but it is not suitable for low Reynolds number flows. \textcolor{black}{We underline that our time discretization combined with our choice of numerical flux in the convection term is not energy conserving. 
} \textcolor{black}{However, we will show in Section \ref{sec:results} that the energy dissipation is remarkably small and it actually helps in avoiding spurious oscillations arising from under resolved scales.}

\section{Computational results}\label{sec:results} 
In this section we validate the proposed method against some well-known benchmark problems. The scheme has been implemented in the finite element library \texttt{NGSolve} \cite{ngSolve}. In the figures, the new scheme is referred to as \lq\lq DG-FEEC\rq\rq, since it is based on both discontinuous Galerkin (DG) methods and compatible finite element exterior calculus (FEEC). \color{black}{Unless otherwise specified, the time-step is computed with the CFL condition 
\begin{equation*}
	\dt = C_{\mathrm{CFL}}\frac{h}{(2r+1)\sigma},
\end{equation*} 
with $C_{\mathrm{CFL}}$ the Courant number, $h$ being a characteristic mesh size and $\sigma\doteq \max(\lVert u\rVert_{\infty}, 1)$.}
If not stated otherwise, $C_{\mathrm{CFL}} = 0.25$ and the physical parameters are $\gamma = 1.4$ and $c_v = 2.5$. We say that we employ polynomials of degree $r$ to as a short-hand for the finite element spaces $\Sigma_{r+1}$, $\RT_r$ and $\DG_r$. \color{black} The symmetric positive definite linear systems are solved with the sparse Cholesky factorization available in \texttt{NGSolve} \cite{ngSolve}.
\color{black}

\subsection{Isentropic vortex}
{\color{black}
To validate the spatial accuracy of the proposed method in the compressible regime, we consider the isentropic vortex proposed by Hu and Shu \cite{HuShuVortex1999}. The domain is the square $\Omega = [0,10]^2$ with periodic boundary conditions. The stationary solution is 
\begin{align*}
\density(\mathbf{x}, t) &= (1 + \delta T)^{\frac{1}{\gamma - 1}},\\
\pres(\mathbf{x}, t) &= (1+ \delta T)^{ \frac{\gamma}{\gamma -1}}, \\
\vel(\mathbf{x}, t) &=  \frac{5}{2\pi}e^{\frac{1-r^2}{2}} (5- y, x -5),
\end{align*}
with $\delta T(\mathbf{x}, t) = -(\gamma-1)\frac{25}{8\gamma\pi^2}e^{1 - r^2}$ and $r = \sqrt{(x-5)^2 + (y-5)^2}$.  We run the simulation until $\tEnd = 1$ with a sequence of unstructured meshes M$_N$ with $N = 40, 60, 80, 100, 120$ being the number of intervals on each side of the square, and we consider polynomial degrees $r = 0,1,2$. The resulting errors and convergence rates are shown in Tables \ref{tab:Shu_errors_r0}, \ref{tab:Shu_errors_r1} and \ref{tab:Shu_errors_r2} respectively. As expected, the convergence rate is $r+1$ when polynomials of degree $r$ are used.
\begin{table}[!hp]
	\caption{Spatial $L^{2}$ error norms and convergence rates at time $t=1$ for the Shu vortex benchmark in 2D with $\mom_h\in\RT_0$.}
	\label{tab:Shu_errors_r0} 	
	\renewcommand{\arraystretch}{1.2}
	\begin{center}
		\color{black}{
		\begin{tabular}{ccccccc}
			\hline 
			Mesh      	&$L^{2}_{\Omega}\left(\densityh \right)$  & $\mathcal{O}\left(\densityh\right)$                       
			&$L^{2}_{\Omega}\left(\velh\right)$  & $\mathcal{O}\left(\velh\right)$&$L^{2}_{\Omega}\left(\presh \right)$ & $\mathcal{O}\left(\presh \right)$ 
			\\ \hline
    M$_{40}$ & $1.1786\cdot 10^{-1}$&$$&$4.0867\cdot 10^{-1}$&$$&$1.4399\cdot 10^{-1}$&$$\\
M$_{60}$ & $8.0714\cdot 10^{-2}$&$0.93$&$2.8752\cdot 10^{-1}$&$0.87$&$9.8315\cdot 10^{-2}$&$0.94$\\
M$_{80}$ & $6.1388\cdot 10^{-2}$&$0.95$&$2.1672\cdot 10^{-1}$&$0.98$&$7.4789\cdot 10^{-2}$&$0.95$\\
M$_{100}$ & $4.9716\cdot 10^{-2}$&$0.95$&$1.7598\cdot 10^{-1}$&$0.93$&$6.0597\cdot 10^{-2}$&$0.94$\\
M$_{120}$ & $4.1559\cdot 10^{-2}$&$0.98$&$1.4757\cdot 10^{-1}$&$0.97$&$5.0599\cdot 10^{-2}$&$0.99$\\
			\hline 
		\end{tabular}}
	\end{center}
\end{table}
\begin{table}[!hp]
\caption{Spatial $L^{2}$ error norms and convergence rates at time $t=1$ for the Shu vortex benchmark in 2D with $\mom_h\in\RT_1$.}
\label{tab:Shu_errors_r1} 	
\renewcommand{\arraystretch}{1.2}
\begin{center}
\color{black}{
\begin{tabular}{ccccccc}
\hline 
Mesh      	&$L^{2}_{\Omega}\left(\densityh \right)$  & $\mathcal{O}\left(\densityh\right)$                       
&$L^{2}_{\Omega}\left(\velh\right)$  & $\mathcal{O}\left(\velh\right)$&$L^{2}_{\Omega}\left(\presh \right)$ & $\mathcal{O}\left(\presh \right)$ 
\\ \hline
M$_{40}$ & $3.3423\cdot 10^{-3}$&$$&$1.0999\cdot 10^{-2}$&$$&$4.2149\cdot 10^{-3}$&$$\\
M$_{60}$ & $1.4795\cdot 10^{-3}$&$2.01$&$4.8212\cdot 10^{-3}$&$2.03$&$1.8802\cdot 10^{-3}$&$1.99$\\
M$_{80}$ & $8.6110\cdot 10^{-4}$&$1.88$&$2.6915\cdot 10^{-3}$&$2.03$&$1.0898\cdot 10^{-3}$&$1.90$\\
M$_{100}$ & $5.4473\cdot 10^{-4}$&$2.05$&$1.7498\cdot 10^{-3}$&$1.93$&$6.9309\cdot 10^{-4}$&$2.03$\\
M$_{120}$ & $3.5057\cdot 10^{-4}$&$2.42$&$1.1496\cdot 10^{-3}$&$2.30$&$4.4690\cdot 10^{-4}$&$2.41$\\
\hline 
\end{tabular}}
\end{center}
\end{table}}
\color{black}
\begin{table}[!hp]
	\caption{Spatial $L^{2}$ error norms and convergence rates at time $t=1$ for the Shu vortex benchmark in 2D with $\mom_h\in\RT_2$.}
	\label{tab:Shu_errors_r2} 	
	\renewcommand{\arraystretch}{1.2}
	\begin{center}
		\color{black}{
		\begin{tabular}{ccccccc}
			\hline 
			Mesh      	&$L^{2}_{\Omega}\left(\densityh \right)$  & $\mathcal{O}\left(\densityh\right)$                       
			&$L^{2}_{\Omega}\left(\velh\right)$  & $\mathcal{O}\left(\velh\right)$&$L^{2}_{\Omega}\left(\presh \right)$ & $\mathcal{O}\left(\presh \right)$ 
			\\ \hline
			M$_{40}$ & $1.6147\cdot 10^{-4}$&$$&$8.2252\cdot 10^{-4}$&$$&$1.9950\cdot 10^{-4}$&$$\\
			M$_{60}$ & $4.5158\cdot10^{-5}$&$3.14$&$2.5168\cdot 10^{-4}$&$2.92$&$5.6302\cdot10^{-5}$&$3.12$\\
			M$_{80}$ & $1.8496\cdot10^{-5}$&$3.10$&$1.0639\cdot 10^{-4}$&$2.99$&$2.3203\cdot10^{-5}$&$3.08$\\
			M$_{100}$ & $9.4169e-06$&$3.03$&$5.4099\cdot10^{-5}$&$3.03$&$1.1721\cdot10^{-5}$&$3.06$\\
			M$_{120}$ & $5.2539e-06$&$3.20$&$2.9573\cdot10^{-5}$&$3.31$&$6.6115\cdot10^{-6}$&$3.14$\\
			\hline 
			\end{tabular}}
			\end{center}
		\end{table}
\subsection{Smooth acoustic wave}
\begin{figure}
\centering
\includegraphics[trim = {5 5 5 5}, clip,width=0.5\textwidth]{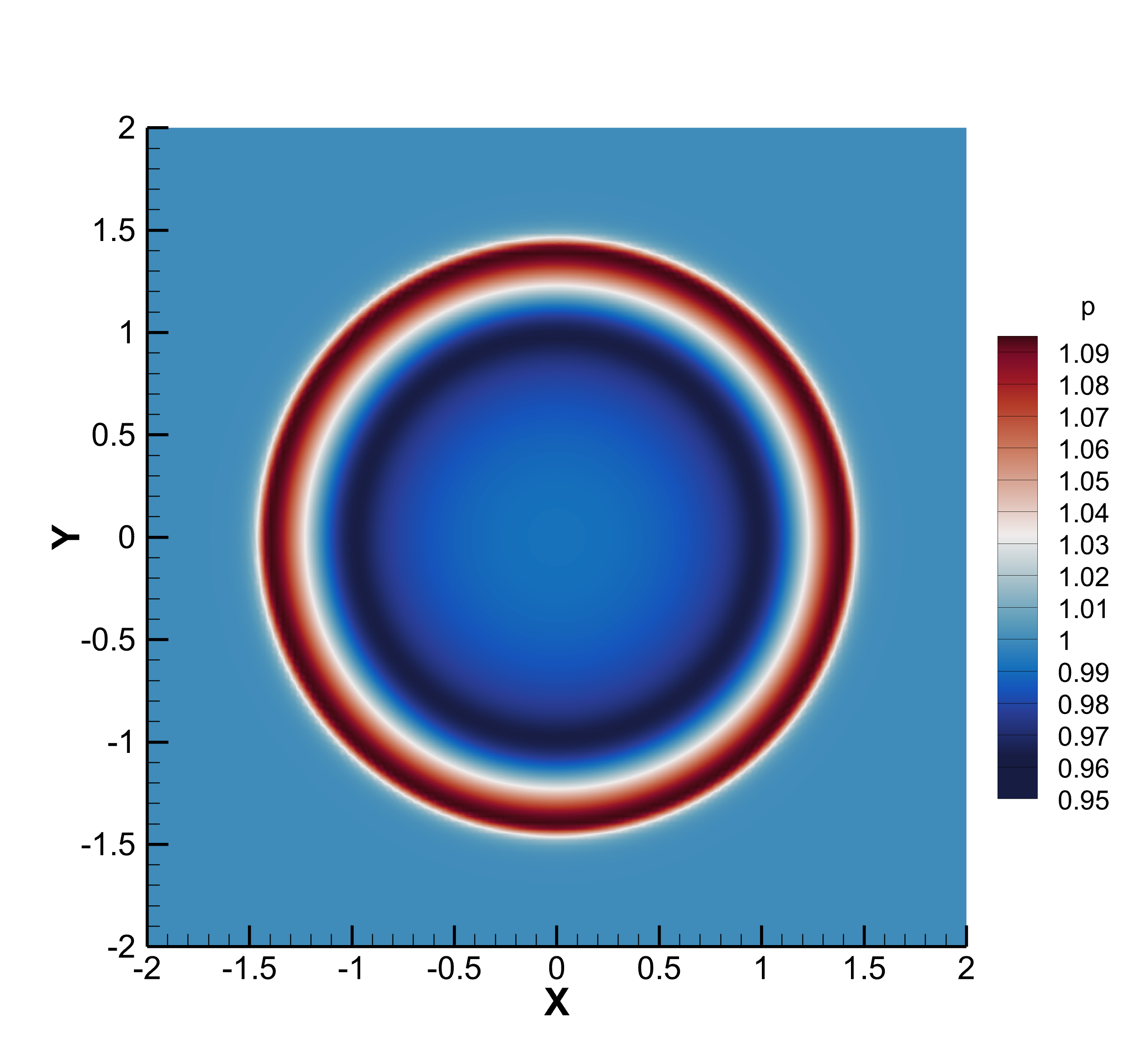}
\caption{Pressure at time $t = 1$ for the smooth acoustic wave test.}
\label{fig:SAW}
\end{figure}

\begin{figure}
\centering
\includegraphics[trim = {5 5 5 5}, clip,width=0.3\textwidth]{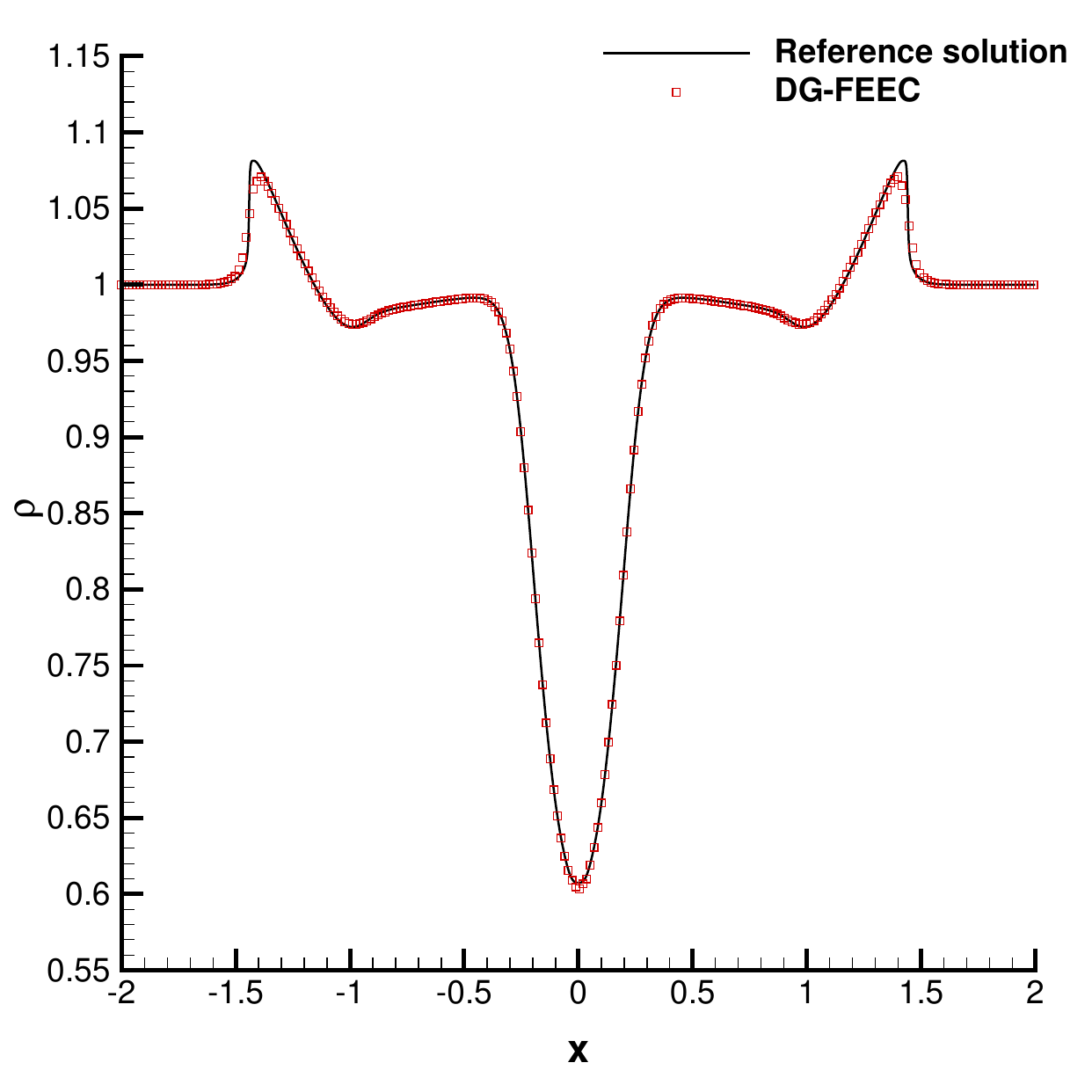}
\includegraphics[trim = {5 5 5 5}, clip,width=0.3\textwidth]{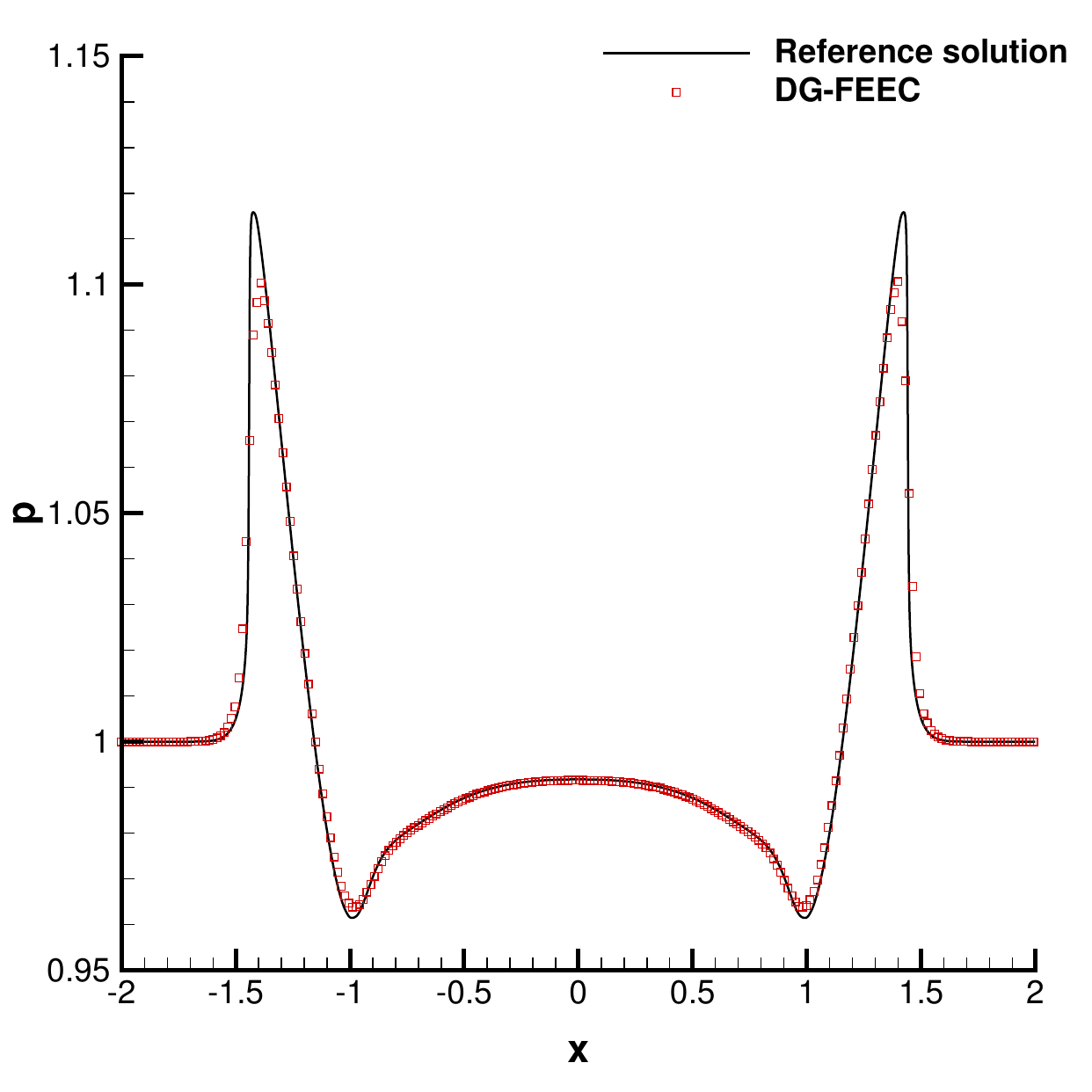}
\includegraphics[trim = {5 5 5 5}, clip,width=0.3\textwidth]{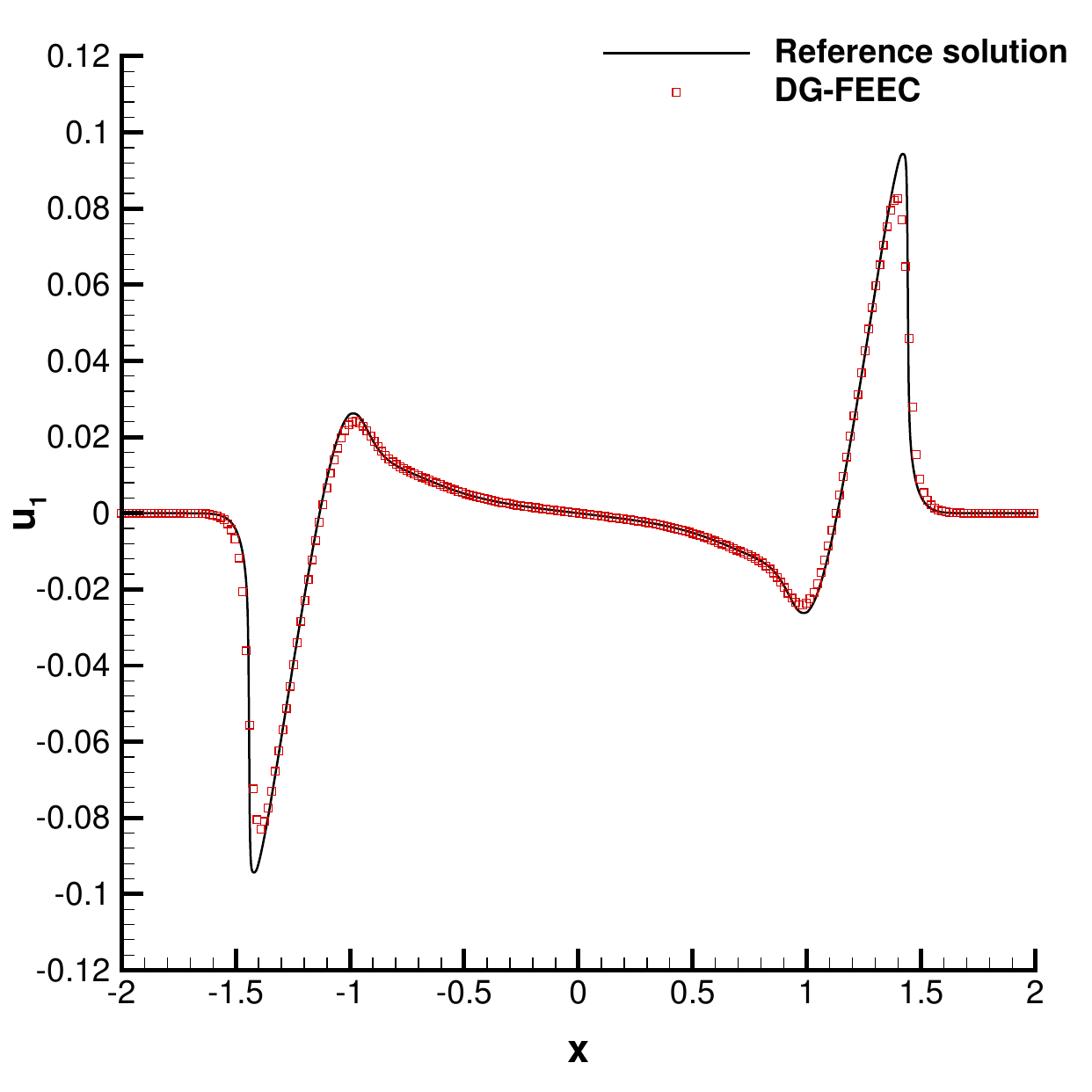}
\caption{Density, pressure and radial velocity along $x = 0$ for the smooth acoustic wave test. Comparison between the DG-FEEC method and the reference solution.}
\label{fig:SAW_cut}
\end{figure}

\color{black}
We consider now the propagation of the wave given at time $t = 0$ by 
\begin{equation*}
\pres(\mathbf{x}, 0) = 1 + e^{-\alpha r^2}, \, \density(\mathbf{x}, 0) = 1, \vel(\mathbf{x}, 0) = \mathbf{0},
\end{equation*}
with $\alpha = 40$ and $r = \sqrt{x^2 + y^2}$. The domain is the periodic square $\Omega = [-2,2]^2$, and is discretized with a $120\times 120$ unstructured triangular mesh. For this test we choose $C_{\mathrm{CFL}}= 0.1$ and we use polynomials of degree $1$. The pressure at the final time $\tEnd = 1$ is plotted in Figure \ref{fig:SAW}. For this test, we can compute a reference solution by solving the equations in the radial direction with a second-order explicit finite volume scheme, see also \cite{TD17} and \cite{Hybrid1}. A comparison with the reference solution and the proposed method with polynomaials of degree $1$ and a $120\times 120$ mesh is shown in Figure \ref{fig:SAW_cut}.
A good agreement is observed between the proposed methodology and the reference solution.

\subsection{Circular explosion}
\begin{figure}
\centering
\includegraphics[trim = {5 5 5 5}, clip,width=0.5\textwidth]{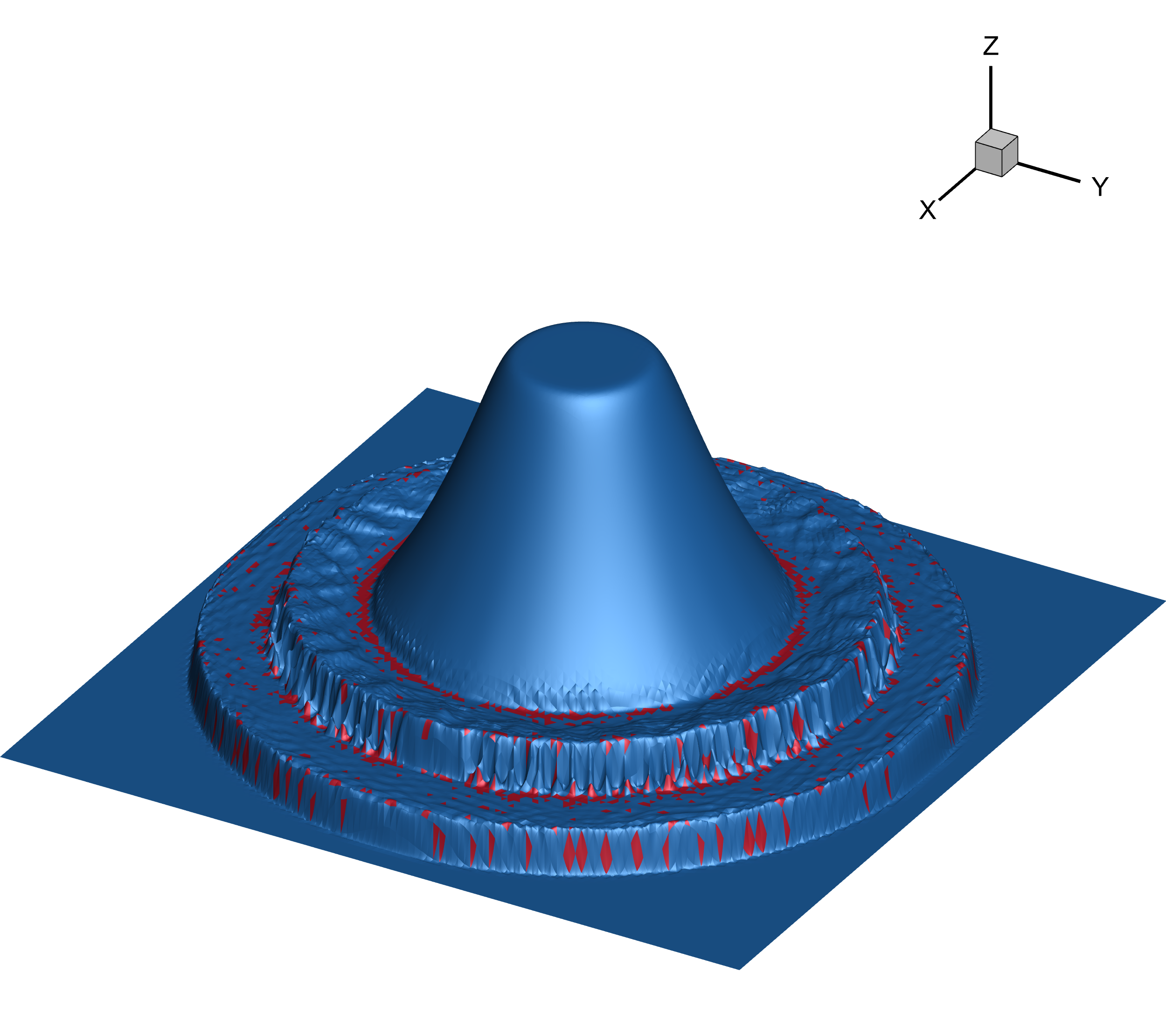}
\caption{Density at time $t = 0.25$ for the circular explosion problem. The artificial viscosity is applied only on the red elements.}
\label{fig:CE}
\end{figure}

\begin{figure}
\centering
\includegraphics[trim = {5 5 5 5}, clip,width=0.3\textwidth]{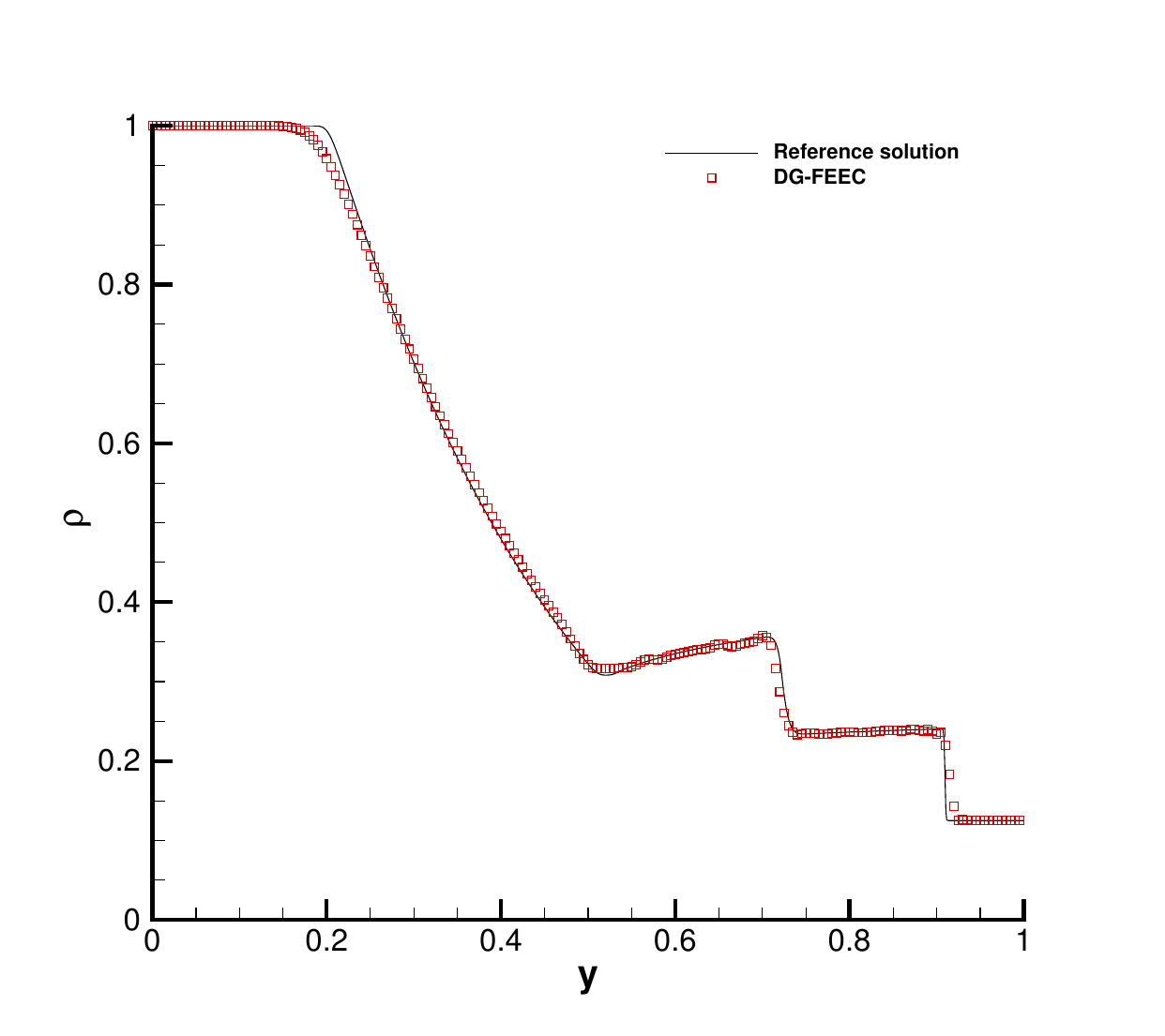}
\includegraphics[trim = {5 5 5 5}, clip,width=0.3\textwidth]{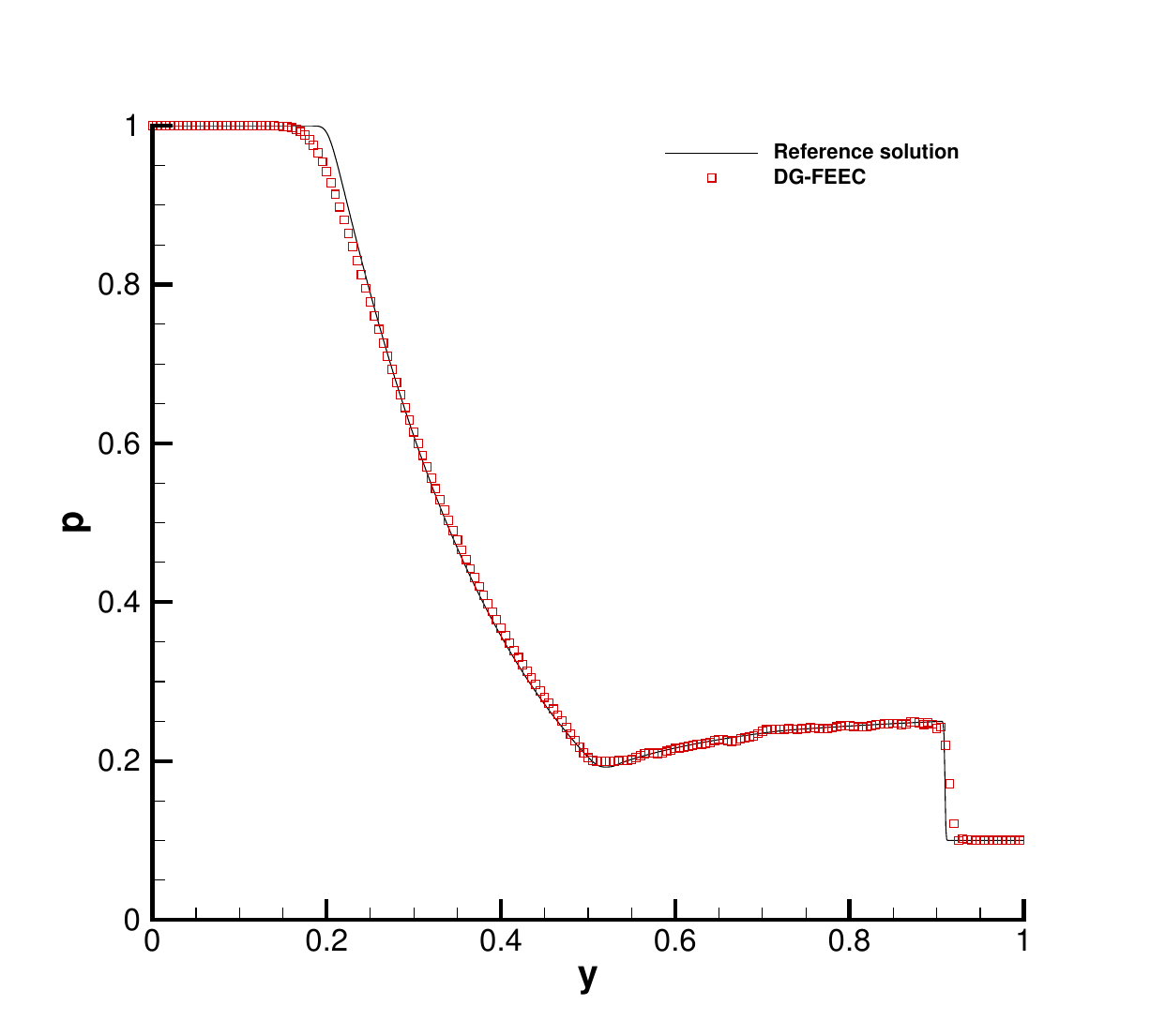}
\includegraphics[trim = {5 5 5 5}, clip,width=0.3\textwidth]{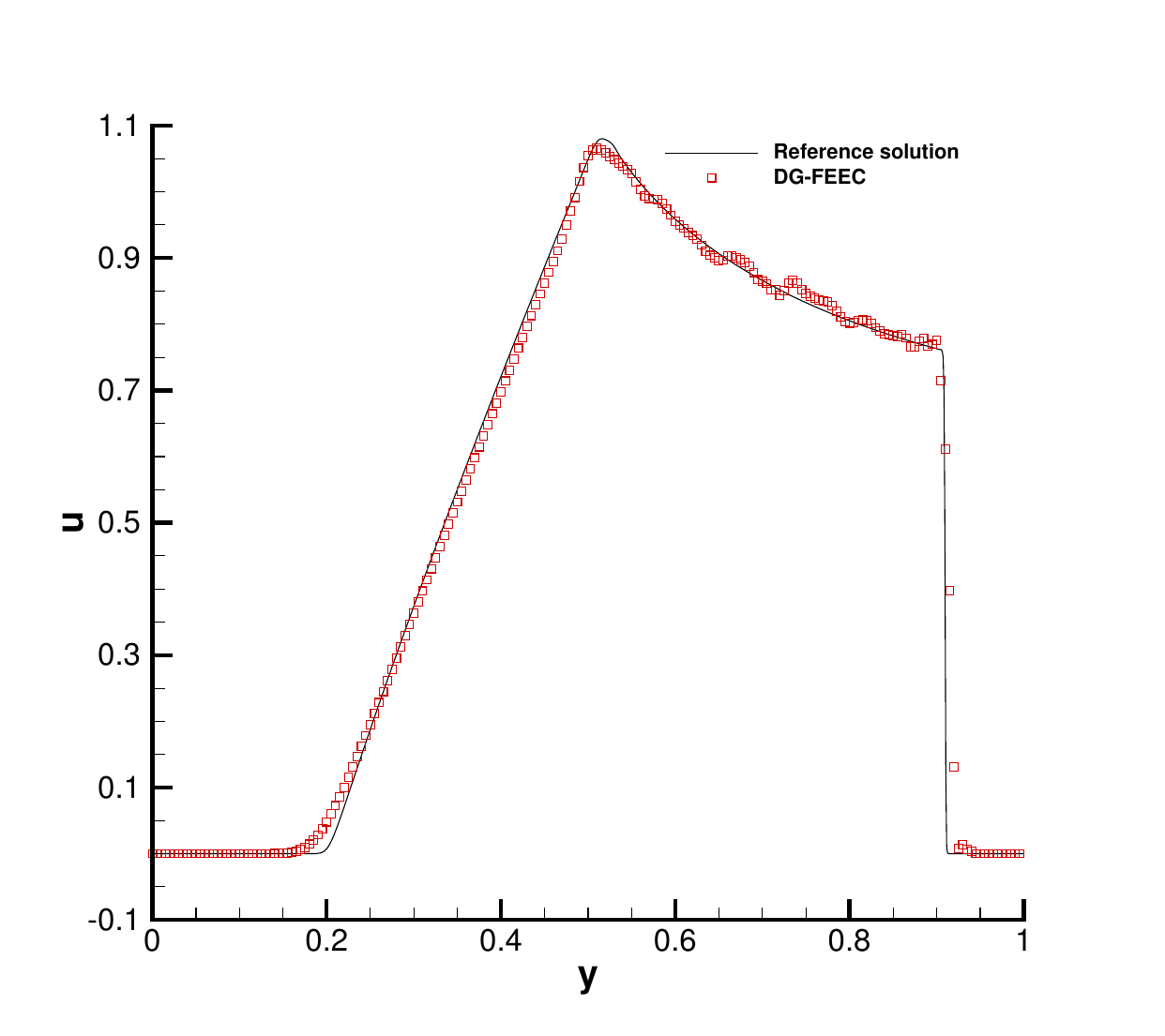}
\caption{Density, pressure and radial velocity along $x = 0$ for the circular explosion test. Comparison between the reference solution and the new DG-FEEC method with polynomial approximation degree two and \textit{a posteriori} artificial viscosity.}
\label{fig:CE_cut}
\end{figure}

\color{black}
To test the robustness of the proposed method on shocks and moderate Mach number flows, we consider now a two-dimensional circular explosion problem (see e.g. \cite{ToroBook, TitarevToro05}). The computational domain is the square $\Omega = [-1,1]^2$ with periodic boundary conditions. The initial data are given by 
\begin{equation*}
\density(\mathbf{x}, 0) = \begin{cases} 1&\text{ if $r\leq 0.5$,}\\0.125&\text{ if $r>0.5$,}\end{cases} \qquad \pres(\mathbf{x}, 0) = \begin{cases} 1&\text{ if $r\leq 0.5$,}\\0.1&\text{ if $r>0.5$,}\end{cases}\qquad \vel(\mathbf{x}, 0) = \mathbf{0}.
\end{equation*}
We run the simulation using polynomials of degree $2$ until $\tEnd = 0.25$ on a $120\times 120$ unstructured mesh. The density at the final time is swhon in Figure \ref{fig:CE}, together with the troubled elements on which the artificial viscosity was applied. As for the smooth acoustic wave test, we have computed a reference solution by solving the equations in the radial direction, 
see \cite{toro-book}. A comparison between the reference solution and the proposed method is shown in Figure \ref{fig:CE_cut}. 

We observe that the a posteriori artificial viscosity is applied only near discontinuities, so that the method does not exibit suprious oscillations while mantaining the high order accuracy. 
\color{black}
\subsection{Kelvin-Helmholtz instability at low Mach}
\begin{figure}
\centering
\includegraphics[trim = {5 5 5 5}, clip,width=0.45\textwidth]{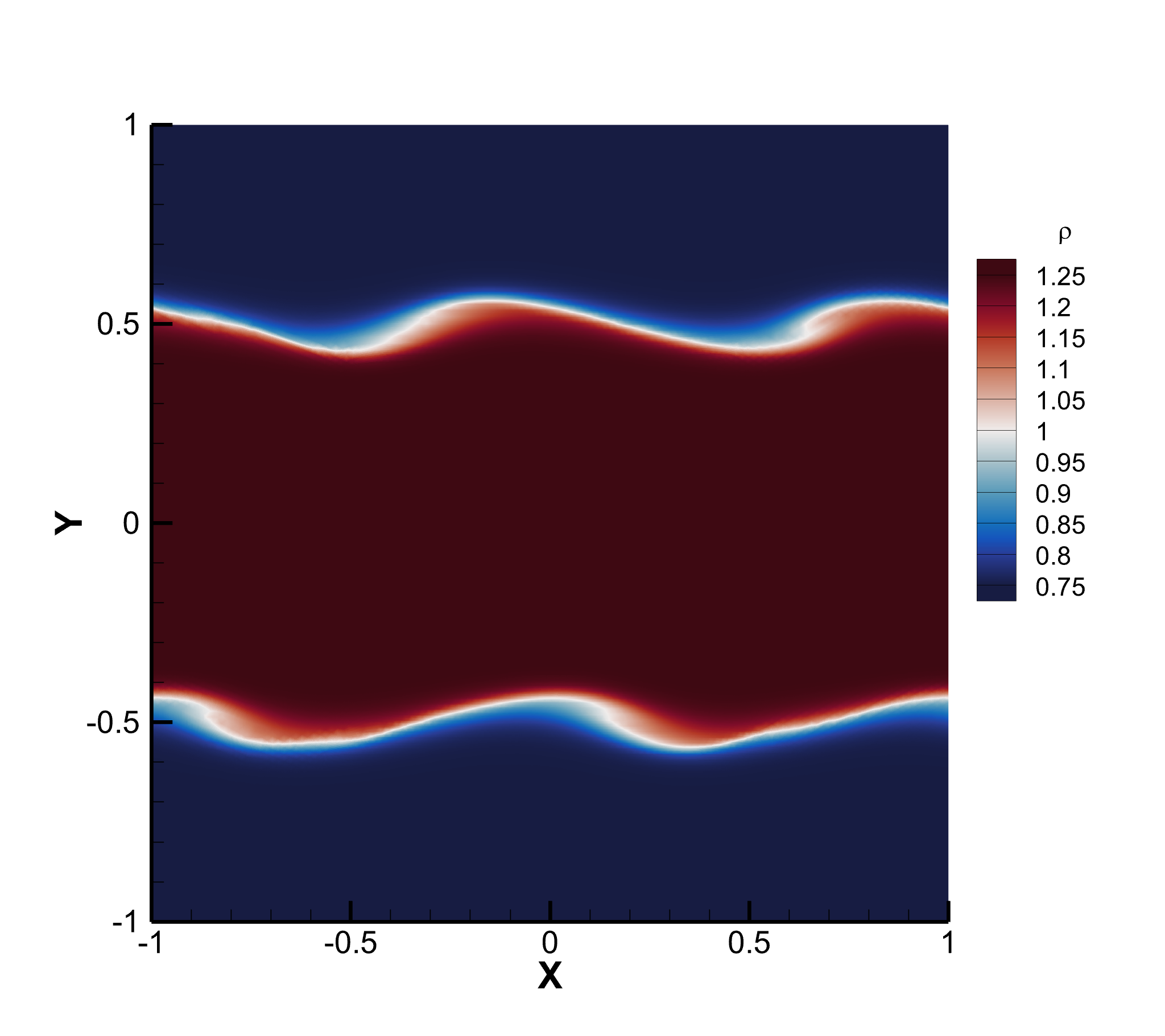}
\includegraphics[trim = {5 5 5 5}, clip,width=0.45\textwidth]{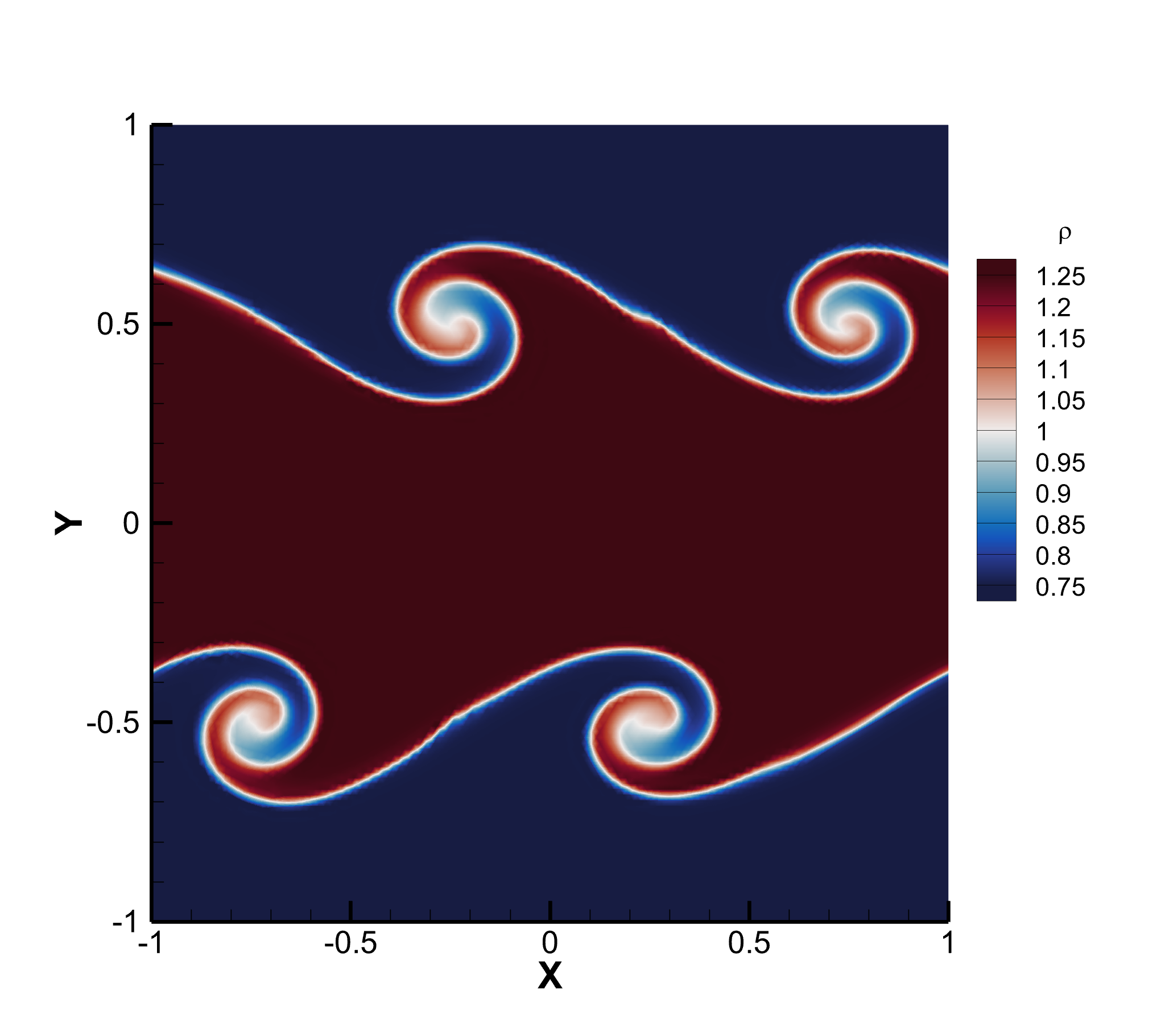}
\includegraphics[trim = {5 5 5 5}, clip,width=0.45\textwidth]{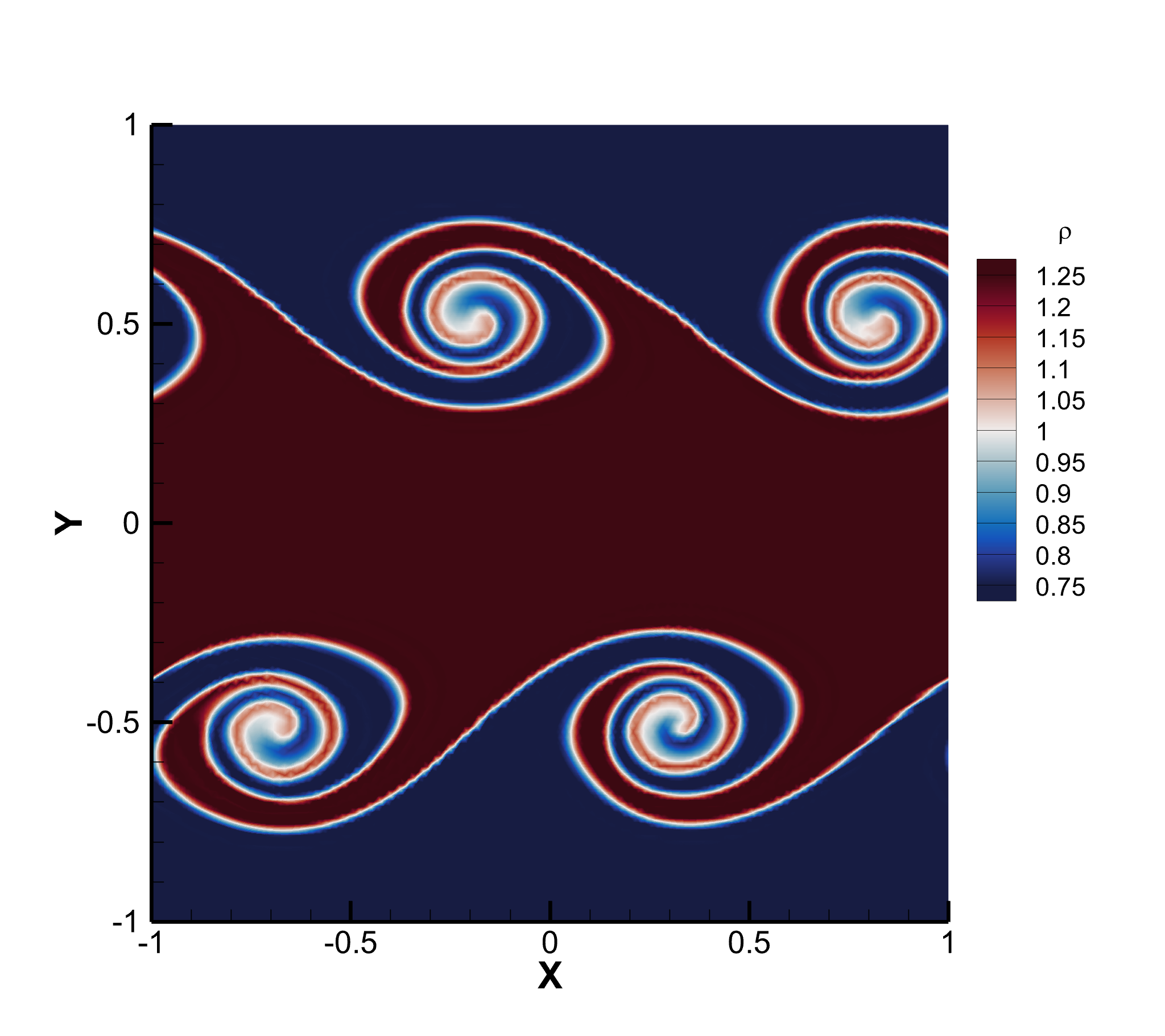}
\includegraphics[trim = {5 5 5 5}, clip,width=0.45\textwidth]{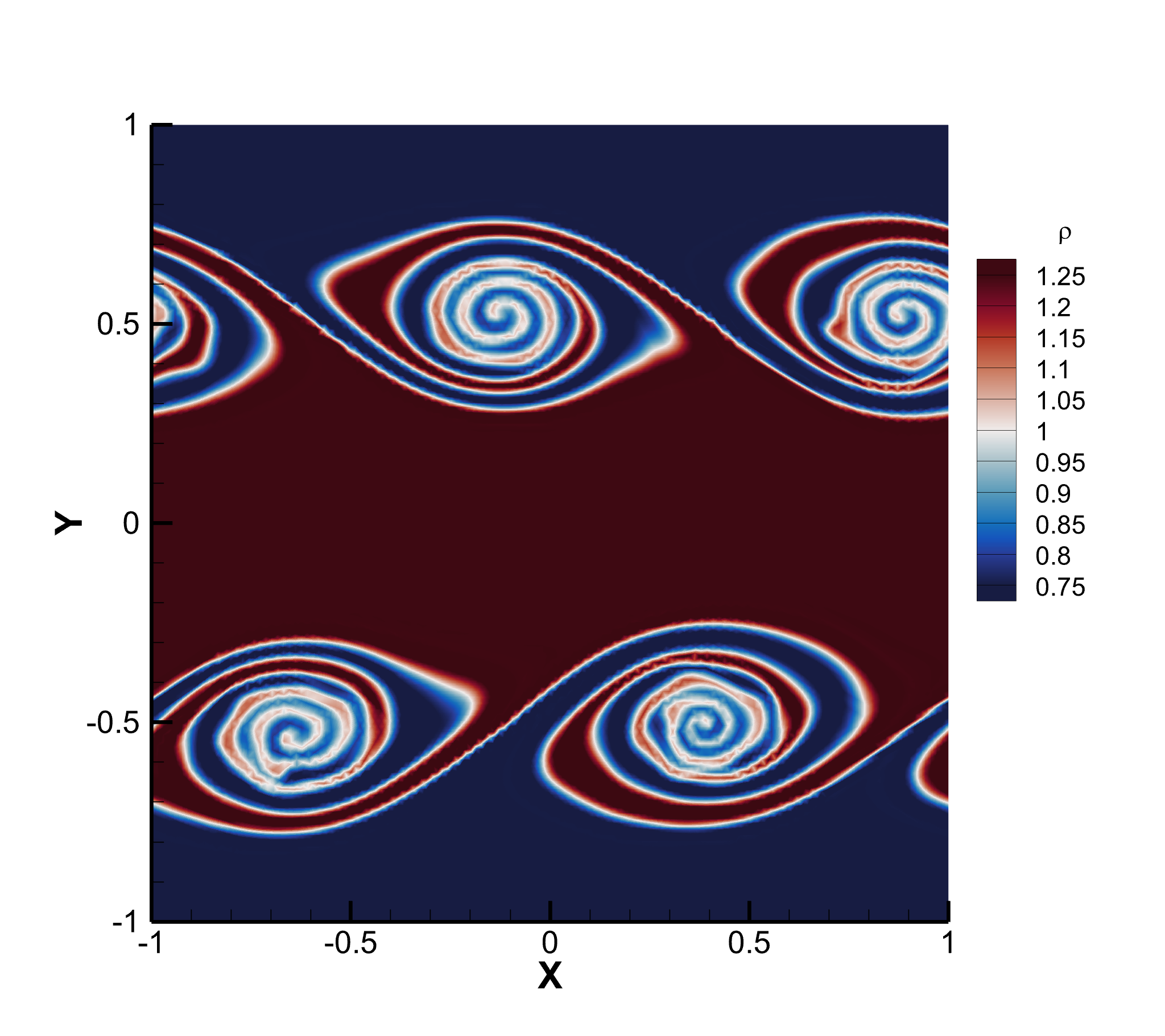}
\caption{Density at times $t = 2,3,4,5$ for the Kelvin-Helmholtz instability test.}
\label{fig:KH}
\end{figure}

We consider now the Kelvin-Helmholtz instability test on the periodic square $\Omega = [-1,1]^2$. The initial condition for this test is 
\begin{alignat*}{2}
\density(\xx,0) &= 1 - \frac{1}{4}\tanh\left(25\left(\lvert y \rvert -\frac{1}{2}\right)\right), \qquad
&&\pres(\xx, 0) = \frac{10^4}{\gamma},\\
\vel_1(\xx, 0) &=-\frac{1}{2}\tanh\left(25\left(\lvert y \rvert -\frac{1}{2}\right)\right),\qquad 
&&\vel_2(\xx, 0) = \frac{1}{100}\sin(2\pi x)cos(2\pi y).\end{alignat*}
The domain is discretized with a $120\times 120$ unstructured triangular mesh, and for this test we emply polynomials of degree $1$. The density at times $t =2,3,4,5$ is displayed in Figure \ref{fig:KH} for a qualitative comparison with other references. 
We remark that this test is outside the theoretical framework of this work, since it is in the incompressible regime, but the density is not constant. 
\color{black}
\subsection{Taylor-Green vortex \textcolor{black}{at low Mach}}
To verify the expected convergence rate of the algorithm \textcolor{black}{and the asymptotic-preserving property}, we consider the stationary Taylor-Green  vortex \cite{chorin1}. The computational domain is the periodic box $\Omega=\left[0,2\pi\right]^2$. The exact solution for this test is 
\begin{equation}
	\density(\xx, t) = 1, \qquad \vel \left(\xx,t\right) = \left( \begin{array}{l} 
		\phantom{-}\sin(x)\cos(y) \\ 
		-\cos(x)\sin(y)\end{array} \right), \qquad 
	\pres\left(\xx,t\right) = p_0 + \frac{1}{4} \left(\cos(2x)+\cos(2y) \right).
	\label{eqn.TGV_ex}
\end{equation}
\color{black}
We run this test until $\tEnd = 0.2$ on $50\times 50$ for $p_0\in\{ 5\cdot 10^3, 5\cdot 10^4, \dots, 5\cdot 10^{12} \}$ to show that $\densityh \to 1$ and $\Div \vel \to 0$ when the Mach number $\Mach$ approaches zero. From \cite{KlaMaj,KlaMaj82,klein,KleinMach,MunzDumbserRoller,PM05} we know that this convergence is quadratic with respect to $\Mach$. We verify this property by reporting the $L^{\infty}$ error of $\density_h$ and $\Div \velh$ as a function of $p_0$ and $\Mach$ in Table \ref{tab:TGV_MachConvergence}. As expected, second order is reached with respect to $\Mach$.
\color{black}
\begin{table}[!hp]
\caption{Spatial $L^{\infty}$ error norms and convergence rates of $\Div \velh$ and $\densityh$ with respect to the Mach number $\Mach$ at time $t=0.2$ for the Taylor-Green vortex benchmark in 2D on a $50\times 50$ mesh with polynomials of degree $1$.}
\label{tab:TGV_MachConvergence} 	
\renewcommand{\arraystretch}{1.2}
\begin{center}
{\color{black}
\begin{tabular}{cccccc}
\hline 
$p_0$ & $\Mach$                           
&$L^{\infty}_{\Omega}\left(\Div\velh\right)$ & $\mathcal{O}\left(\Div\velh\right)$ 
&$L^{\infty}_{\Omega}\left(\rho_h \right)$ & $\mathcal{O}\left(\presh \right)$ 
\\ \hline
   $5\cdot 10^{3}$ & $1.20\cdot 10^{-2}$ & $5.4721\cdot 10^{-5}$&$$&$1.0844\cdot 10^{-5}$&$$\\
  $5\cdot 10^{4}$ & $3.78\cdot 10^{-3}$ & $5.4720\cdot 10^{-6}$&$2.00$&$1.0844\cdot 10^{-6}$&$2.00$\\
  $5\cdot 10^{5}$ & $1.20\cdot 10^{-3}$ & $5.4720\cdot 10^{-7}$&$2.00$&$1.0844\cdot 10^{-7}$&$2.00$\\
  $5\cdot 10^{6}$ & $3.78\cdot 10^{-4}$ & $5.4720\cdot 10^{-8}$&$2.00$&$1.0844\cdot 10^{-8}$&$2.00$\\
  $5\cdot 10^{7}$ & $1.20\cdot 10^{-4}$ & $5.4719\cdot 10^{-9}$&$2.00$&$1.0844\cdot 10^{-9}$&$2.00$\\
  $5\cdot 10^{9}$ & $3.78\cdot 10^{-5}$ & $5.4726\cdot 10^{-10}$&$2.00$&$1.0845\cdot 10^{-10}$&$2.00$\\
  $5\cdot 10^{9}$ & $1.20\cdot 10^{-5}$ & $5.4810\cdot 10^{-11}$&$2.00$&$1.0844\cdot 10^{-11}$&$2.00$\\
  $5\cdot 10^{10}$ & $3.78\cdot 10^{-6}$ & $5.4292\cdot 10^{-12}$&$2.01$&$1.0894\cdot 10^{-12}$&$2.00$\\
  $5\cdot 10^{11}$ & $1.20\cdot 10^{-6}$ & $6.3112\cdot 10^{-13}$&$1.87$&$1.0703\cdot 10^{-13}$&$2.02$\\
  $5\cdot 10^{12}$ & $3.78\cdot 10^{-7}$ & $2.1806\cdot 10^{-13}$&$0.92$&$6.6613\cdot 10^{-15}$&$2.41$\\
\hline 
\end{tabular}}
\end{center}
\end{table}
\color{black}
 \textcolor{black}{We now run the simulation until $\tEnd = 0.5$ with the same sequence of meshes used for the Hu-Shu vortex and $p_0 = 10^7$.} We use \textcolor{black}{ polynomials of degree $r$ with $r = 0,1,2$} and we set $C_{\mathrm{CFL}}= 0.5$. The results are shown in Tables~\ref{tab:TGV_errors_r1}, \ref{tab:TGV_errors_r2} and \ref{tab:TGV_errors_r3} respectively. As expected, we observe $r+1$-{th} order accuracy for both velocity and pressure when using \textcolor{black}{polynomials of degree $r$.}

\color{black}

\color{black}

\begin{figure}
	\caption{Evolution of energy (left), momentum (center), incompressibility (right) for the inviscid Taylor-Green vortex for the time interval $t\in [0,10]$. For this test we have used a $40\times 40$ mesh and polynomials of degree $2$.}
	\label{fig:cons} 
	\centering
	\includegraphics[trim = {5 5 5 5}, clip, width =0.3\textwidth]{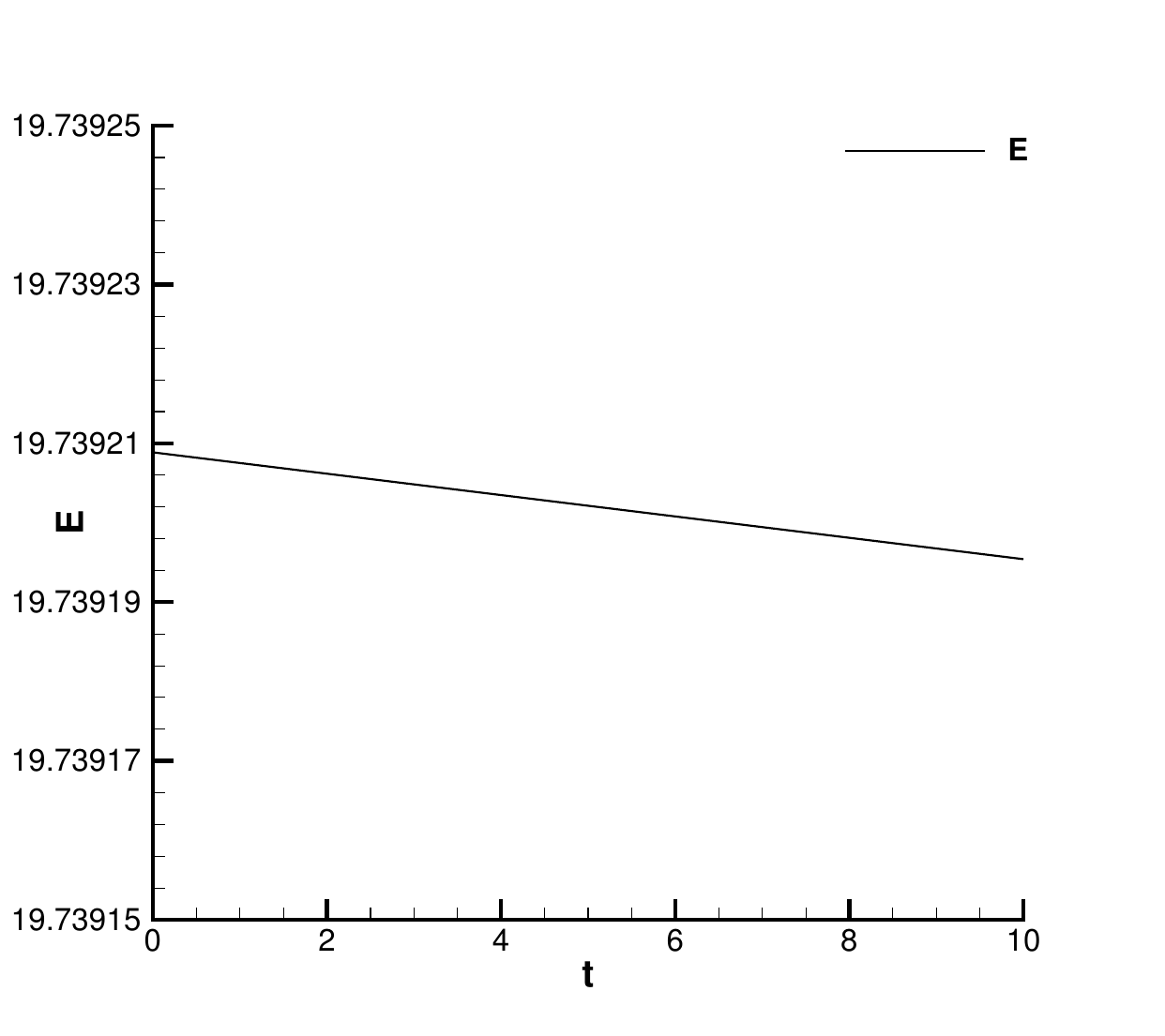}
	\includegraphics[trim = {5 5 5 5}, clip,width =0.3\textwidth]{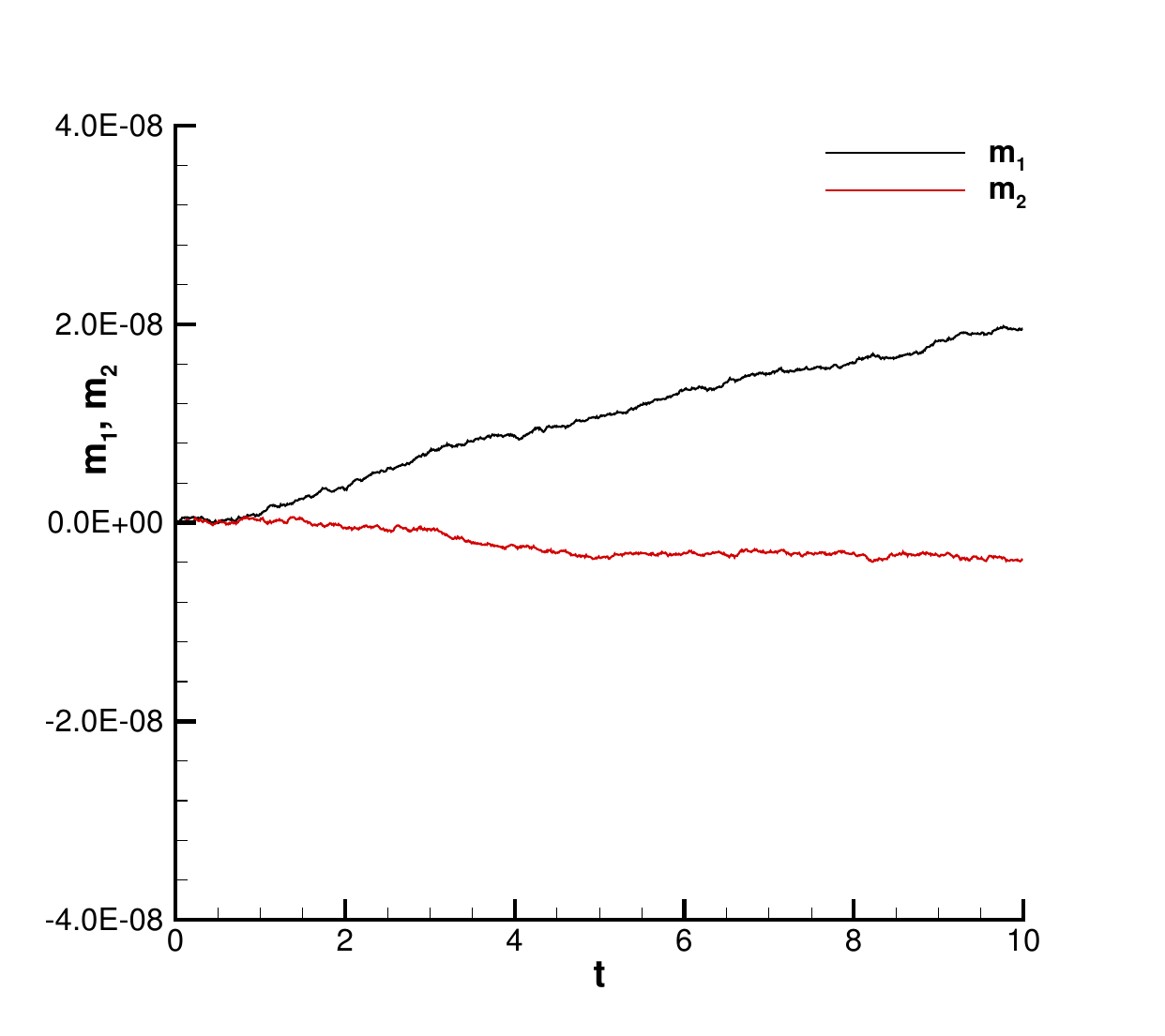}
	\includegraphics[trim = {5 5 5 5}, clip,width =0.3\textwidth]{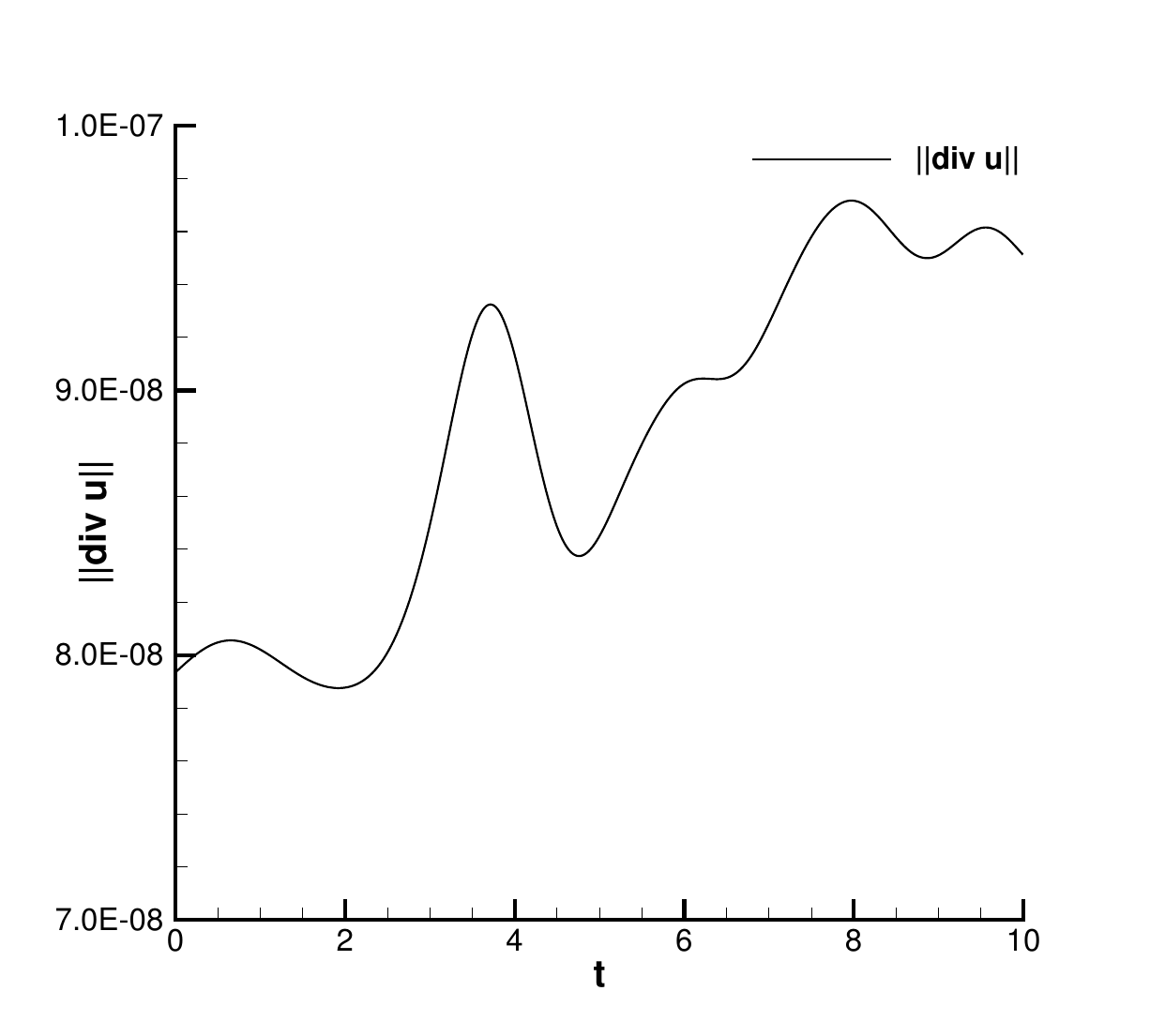}
\end{figure}

 \textcolor{black}{To test the conservation properties of our scheme, we now repeat the test with $p_0 = 1e-7$ until $\tEnd = 10$ using the coarser mesh, $r = 2$, keeping track of total energy $E = \frac{1}{2}\int_{\Omega}\rho\lvert\vel\rvert^2 \dx$, total momentum $m_i = \int_{\Omega}\mom_i\dx$ for $i = 1,2$ and incompressibility $\lVert \Div \vel\rVert$. The evolution of these quantities is displayed in Figure~\ref{fig:cons}. For the energy, we observe a remarkably slow  dissipation rate (less than $\approx 2\cdot 10^{-6}$ per time unit), while momentum and incompressilibity have an error of the same order of magnitude of the square of the Mach number.} 

\color{black}
 
\begin{table}[!hp]
	\caption{Spatial $L^{2}$ error norms and convergence rates at time $t=0.5$ for the Taylor-Green vortex benchmark in 2D with polynomials of degree $0$.}
	\label{tab:TGV_errors_r1} 	
	\renewcommand{\arraystretch}{1.2}
	\begin{center}
	\color{black}
		\begin{tabular}{cccccc}
			\hline 
			Mesh                            
			&$L^{2}_{\Omega}\left(\velh\right)$ & $\mathcal{O}\left(\velh\right)$ 
			&$L^{2}_{\Omega}\left(\presh \right)$ & $\mathcal{O}\left(\presh \right)$ 
			\\ \hline
     M$_{40}$ & $3.6930\cdot 10^{-1}$&$$&$2.2559\cdot 10^{-1}$&$$\\
     M$_{60}$ & $2.4669\cdot 10^{-1}$&$1.00$&$1.5437\cdot 10^{-1}$&$0.94$\\
     M$_{80}$ & $1.8692\cdot 10^{-1}$&$0.96$&$1.1682\cdot 10^{-1}$&$0.97$\\
     M$_{100}$ & $1.5159\cdot 10^{-1}$&$0.94$&$9.2404\cdot 10^{-2}$&$1.05$\\
     M$_{120}$ & $1.2639\cdot 10^{-1}$&$1.00$&$7.8722\cdot 10^{-2}$&$0.88$\\
			\hline 
		\end{tabular}
		\color{black}
	\end{center}
\end{table}

\begin{table}[!hp]
\caption{Spatial $L^{2}$ error norms and convergence rates at time $t=0.5$ for the Taylor-Green vortex benchmark in 2D with polynomials of degree $1$.}
\label{tab:TGV_errors_r2} 	
\renewcommand{\arraystretch}{1.2}
\begin{center}
\color{black}
	\begin{tabular}{ccccccc}
		\hline 
		Mesh                            
		&$L^{2}_{\Omega}\left(\velh\right)$ & $\mathcal{O}\left(\velh\right)$ 
		&$L^{2}_{\Omega}\left(\presh \right)$ & $\mathcal{O}\left(\presh \right)$ 
		\\ \hline
M$_{40}$ & $7.2829\cdot 10^{-3}$&$$&$5.2217\cdot 10^{-3}$&$$\\
M$_{60}$ & $3.2420\cdot 10^{-3}$&$2.00$&$2.2959\cdot 10^{-3}$&$2.03$\\
M$_{80}$ & $1.8173\cdot 10^{-3}$&$2.01$&$1.3107\cdot 10^{-3}$&$1.95$\\
M$_{100}$ & $1.1786\cdot 10^{-3}$&$1.94$&$8.4328\cdot 10^{-4}$&$1.98$\\
M$_{120}$ & $8.2138\cdot 10^{-4}$&$1.98$&$5.7835\cdot 10^{-4}$&$2.07$\\
		\hline 
	\end{tabular}
	\color{black}
\end{center}
\end{table}

\begin{table}[!hp]
	\caption{Spatial $L^{2}$ error norms and convergence rates at time $t=0.5$ for the Taylor-Green vortex benchmark in 2D with polynomials of degree $2$.}
\label{tab:TGV_errors_r3} 	
\renewcommand{\arraystretch}{1.2}
\begin{center}
\color{black}
	\begin{tabular}{ccccccc}
		\hline 
		Mesh                            
		&$L^{2}_{\Omega}\left(\velh\right)$  & $\mathcal{O}\left(\velh\right)$&$L^{2}_{\Omega}\left(\presh \right)$ & $\mathcal{O}\left(\presh \right)$ 
		\\ \hline
   M$_{40}$ & $2.0769\cdot 10^{-4}$&$$&$9.2427\cdot 10^{-5}$&$$\\
M$_{60}$ & $6.1905\cdot 10^{-5}$&$2.99$&$2.6740\cdot 10^{-5}$&$3.06$\\
M$_{80}$ & $2.6161\cdot 10^{-5}$&$2.99$&$1.1173\cdot 10^{-5}$&$3.03$\\
M$_{100}$ & $1.3545\cdot 10^{-5}$&$2.95$&$5.8688\cdot 10^{-6}$&$2.89$\\
M$_{120}$ & $7.8322\cdot 10^{-6}$&$3.00$&$3.4815\cdot 10^{-6}$&$2.86$\\
		\hline 
	\end{tabular}
	\color{black}
\end{center}
\end{table}

\subsection{Double Shear Layer \textcolor{black}{at low Mach}}
\begin{figure}
	\centering
	\includegraphics[trim = {5 5 5 5}, clip,width=0.45\textwidth]{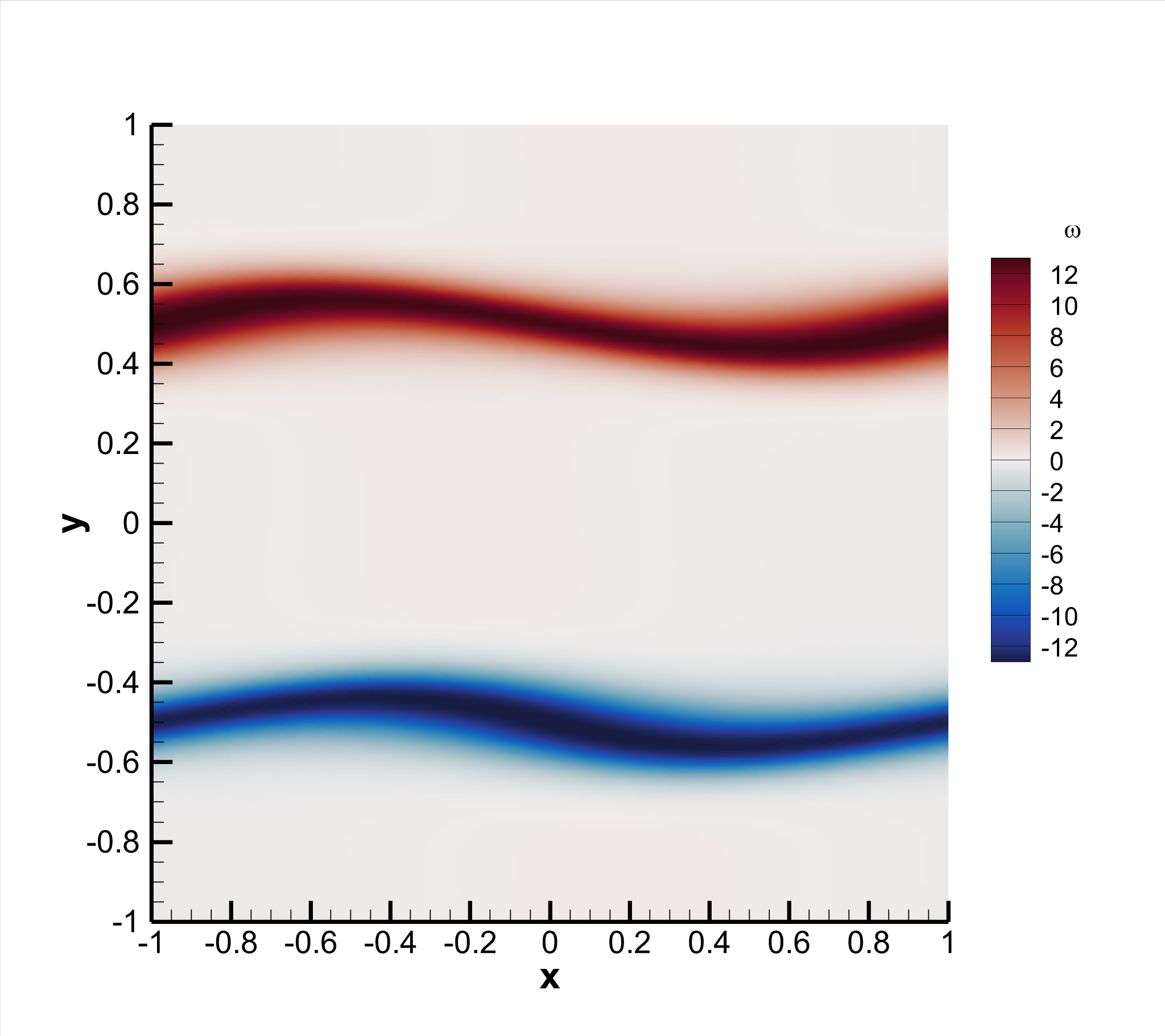}
		\includegraphics[trim = {5 5 5 5}, clip,width=0.45\textwidth]{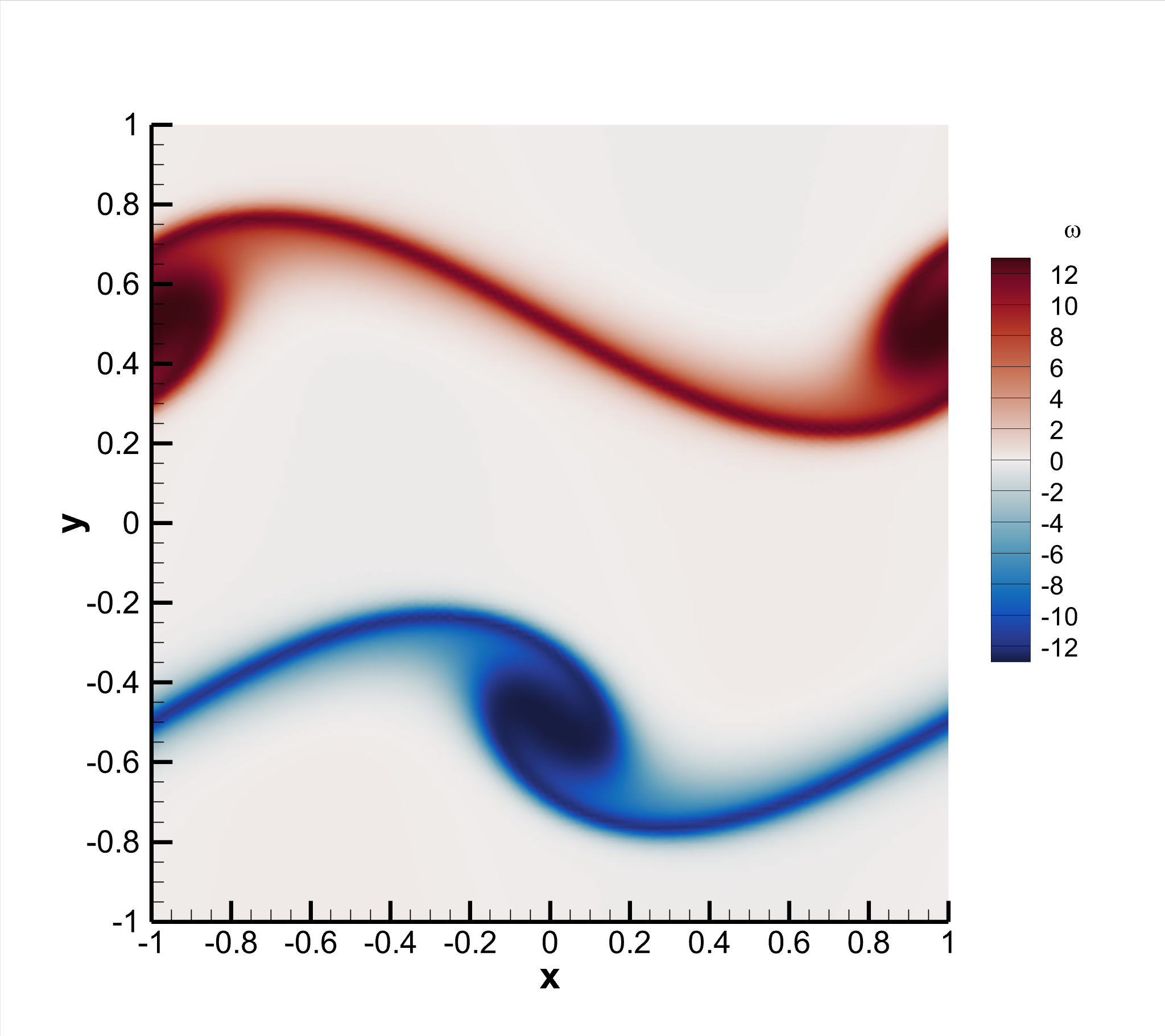}
			\includegraphics[trim = {5 5 5 5}, clip,width=0.45\textwidth]{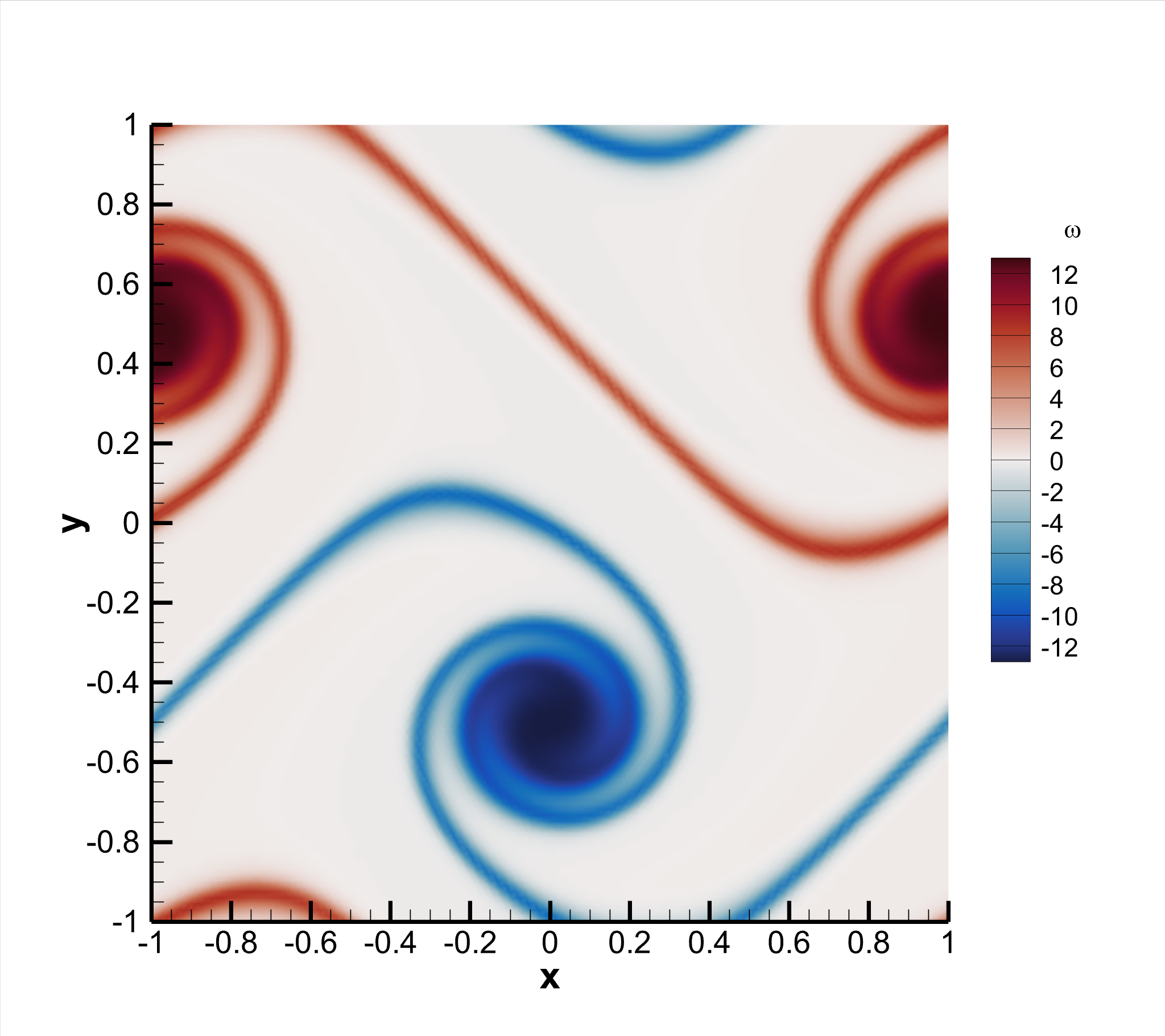}
				\includegraphics[trim = {5 5 5 5}, clip,width=0.45\textwidth]{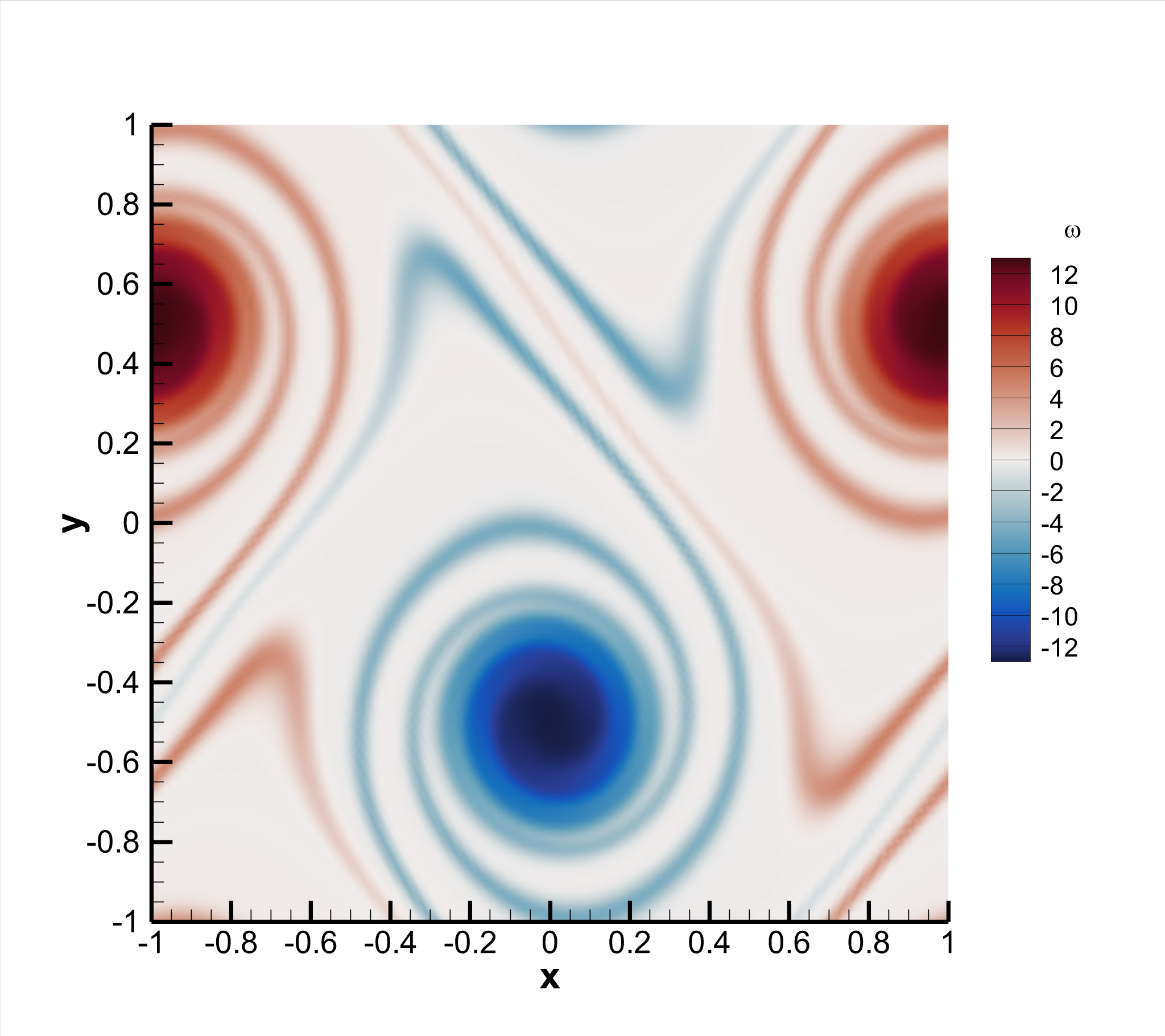}
				\caption{Vorticity at times $t = 0.8$, $1.6$, $2.4$ and $3.6$ for the Double Shear Layer test.}
				\label{fig:DSL}
\end{figure}

We consider now the double shear layer test \cite{BCG89}. For this test the computational domain is $\Omega = [-1,1]^2$ with periodic boundary conditions. The viscosity is set to $\mu = 2\times 10^{-4}$. We consider an initial condition given by
\begin{gather*}
	\density(\xx, t) = 1, \qquad\pres \left(\xx,0\right) = \frac{10^{4}}{\gamma}, \qquad
	\vel_{1} \left(\xx,0\right) = \left\lbrace \begin{array}{lc}
		\tanh \left[\hat{\rho}(\hat{y}-0.25)\right] & \mathrm{ if } \; \hat{y} \le 0.5,\\
		\tanh \left[\hat{\rho}(0.75-\hat{y})\right] & \mathrm{ if } \; \hat{y} > 0.5,
	\end{array}\right.
	\\ 
	\vel_{2} \left(\xx,0\right) = \delta \sin \left(2\pi \hat{x}\right), \qquad
	\hat{x}= \frac{x+1}{2}, \qquad \hat{y}= \frac{y+1}{2},
\end{gather*}
with $\hat{\rho} = 30$ and $\delta =0.05$ being the parameters that determine the slope of the shear layer and the amplitude of the initial perturbation. For this test we use a mesh with $120$ elements on each side and we use a fixed time-step $\dt = 10^{-4}$. The contours of the vorticity $\vorh$ at times $t = 0.8$, $1.6$, $2.4$ and $3.6$ are shown in Figure~\ref{fig:DSL} for a qualitative comparison with other references, see e.g. \cite{BCG89,TD15,Hybrid1,HTCGPR,HTCAbgrall}. 

\subsection{Lid-driven cavity \textcolor{black}{at low Mach}}
\begin{figure}
	\centering
	\includegraphics[trim = {5 5 5 5}, clip, width = 0.45\textwidth]{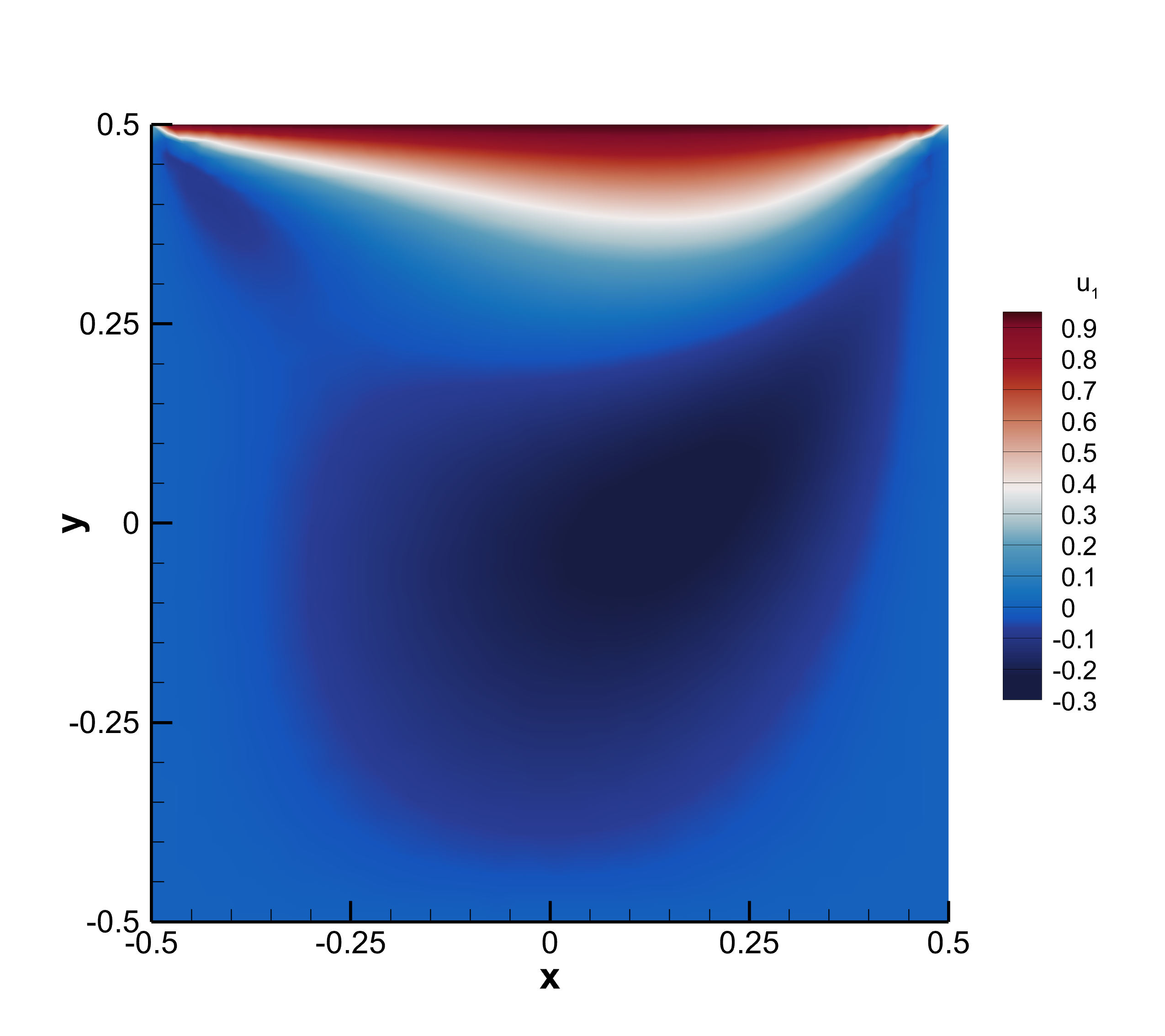}
	\includegraphics[trim = {5 5 5 5}, clip,width = 0.45\textwidth]{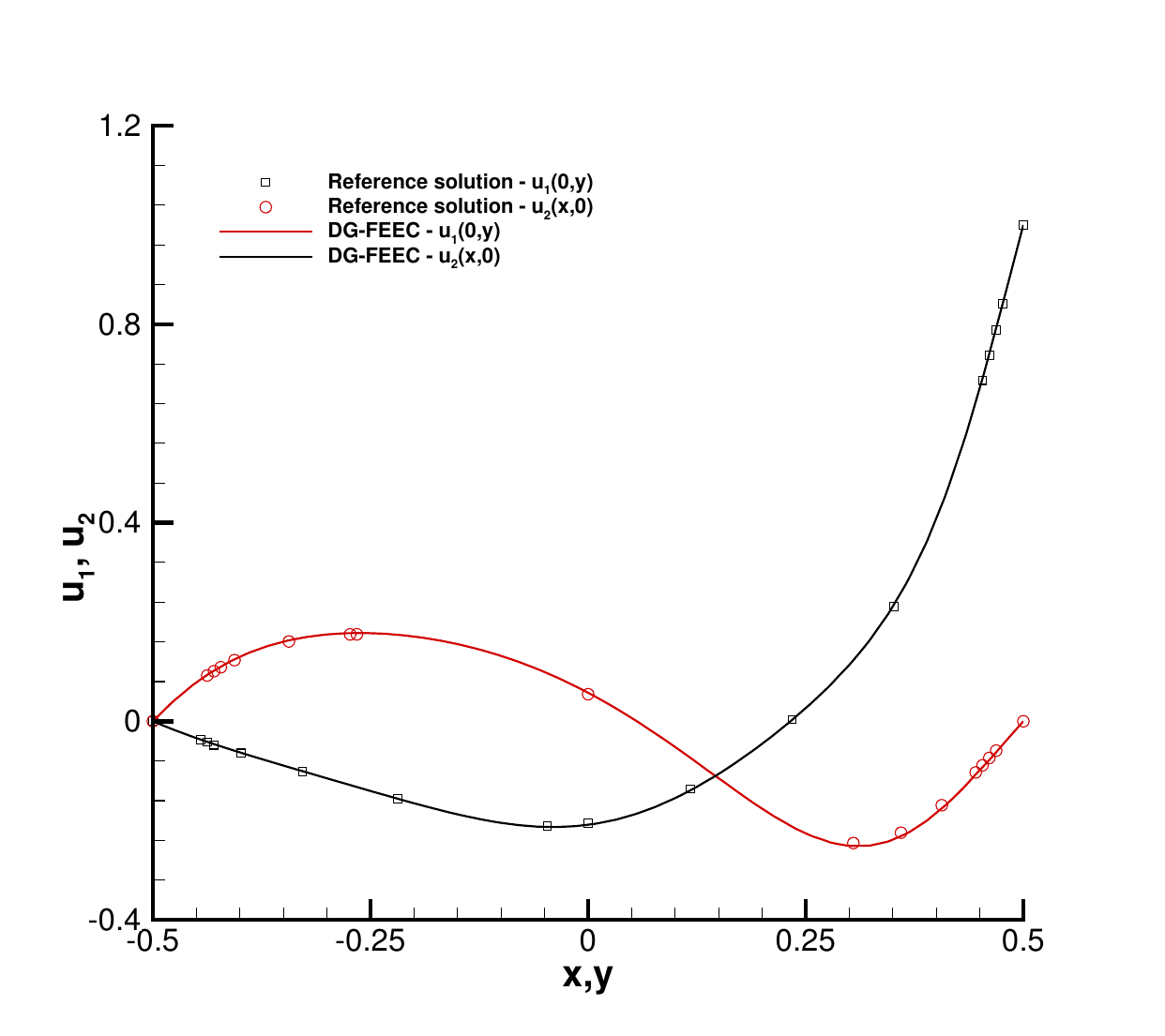}
	\caption{Contour plot of $\vel_1$ (left) for the lid-driven cavity test with $\vis = 0.01$ and comparison with the reference solution of \cite{GGS82} for $\vel_{1}(0,y)$ and $\vel_{2}(x,0)$ (right).}
	\label{fig:LDC}
\end{figure}

\begin{figure}
\centering
\includegraphics[trim = {5 5 5 5}, clip, width = 0.45\textwidth]{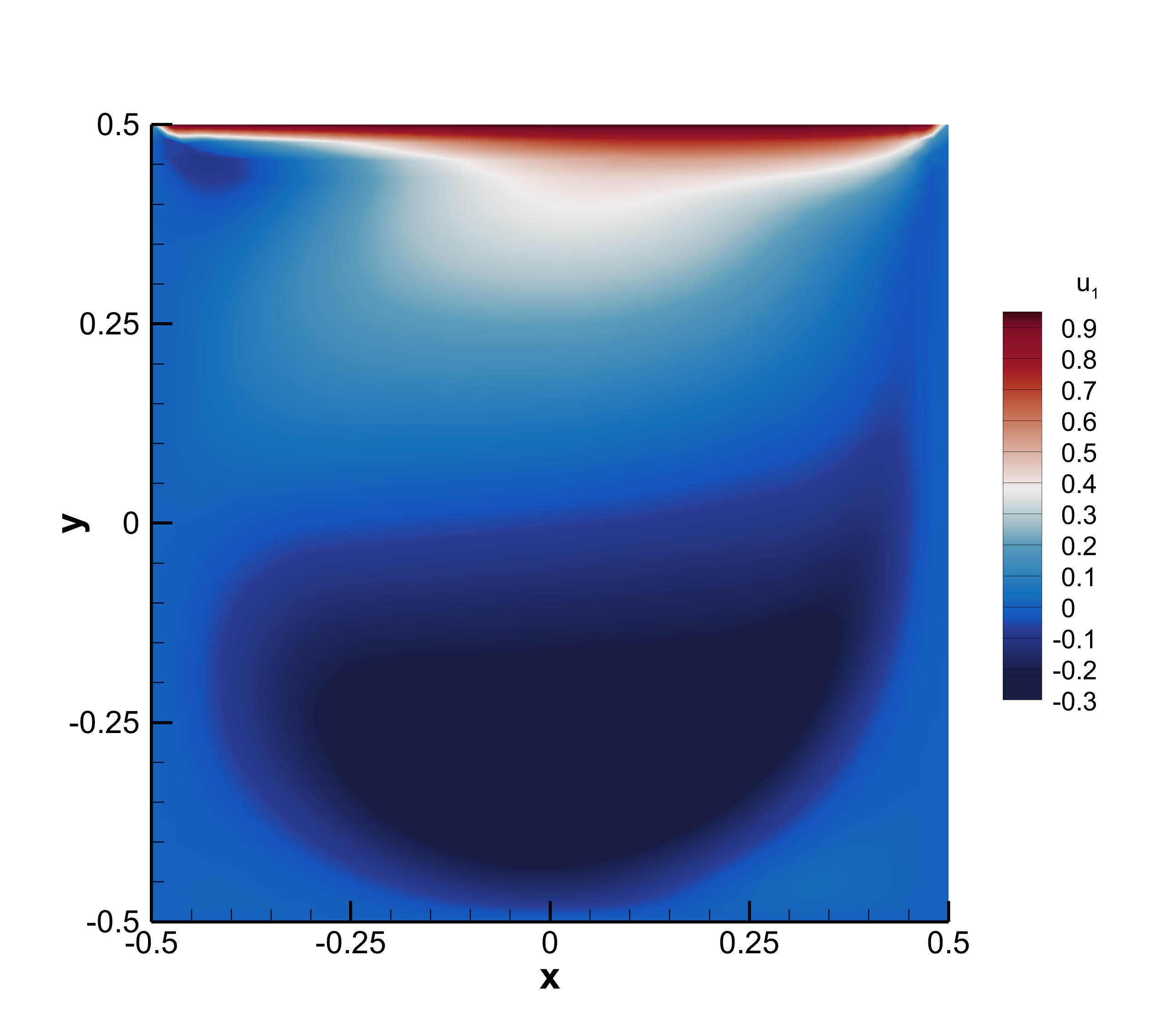}
\includegraphics[trim = {5 5 5 5}, clip,width = 0.45\textwidth]{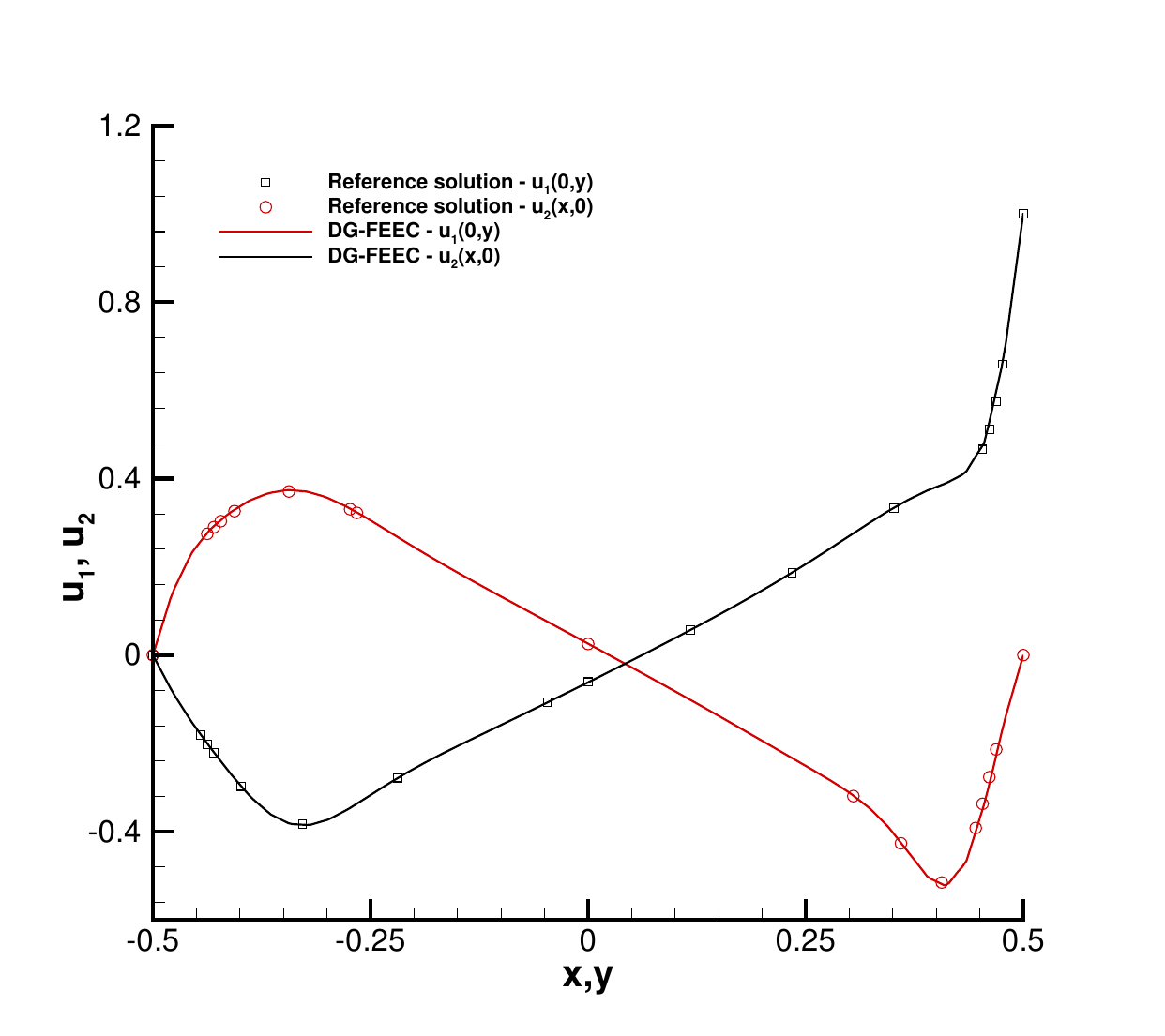}
\caption{Contour plot of $\vel_1$ (left) for the lid-driven cavity test with $\vis = 0.001$ and comparison with the reference solution of \cite{GGS82} for $\vel_{1}(0,y)$ and $\vel_{2}(x,0)$ (right).}
\label{fig:LDC_Re1000}
\end{figure}

We consider the classical lid-driven cavity problem proposed by Ghia, Ghia and Shin \cite{GGS82}. The domain is $\Omega = [-0.5,0.5]^2$ with Dirichlet boundary conditions for the velocity. In particular, we set $\overline{\vel} = (1,0)$ if $y = 0.5$ and $\overline{\vel}= (0,0)$ otherwise. \textcolor{black}{For this test, we set $\mu = 0.01$.} We discretize the domain with a mesh that has $40$ elements on each side and we run the simulation until $\tEnd = 10$. The result is shown in Figure~\ref{fig:LDC}. An excellent agreement with the reference solution \textcolor{black}{from \cite{GGS82}} is observed. 
\color{black} We repeat now the same test with $\mu = 0.001$ and the same mesh, but this time we run the simulation until a steady state is reached. The result is shown in Figure~\ref{fig:LDC_Re1000}. Again, good agreement is observed between our method and the reference solution from \cite{GGS82}.
\color{black}

\subsection{Backward-facing step \textcolor{black}{at low Mach}}
\begin{figure}
	\centering
	\includegraphics[trim = {5 5 5 5}, clip, width = 0.45\textwidth]{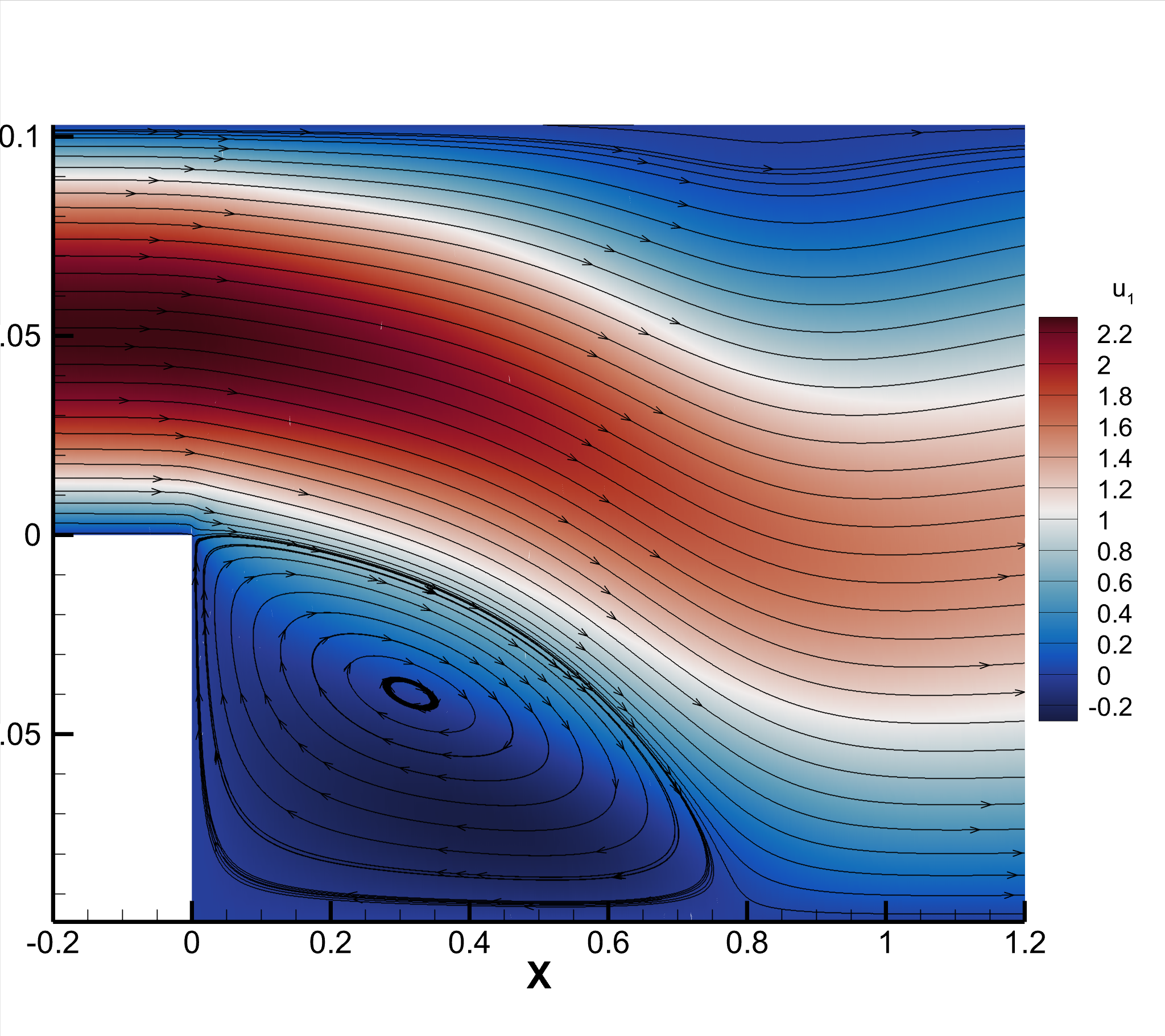}
	\includegraphics[trim = {5 5 5 5}, clip, width=0.45\textwidth]{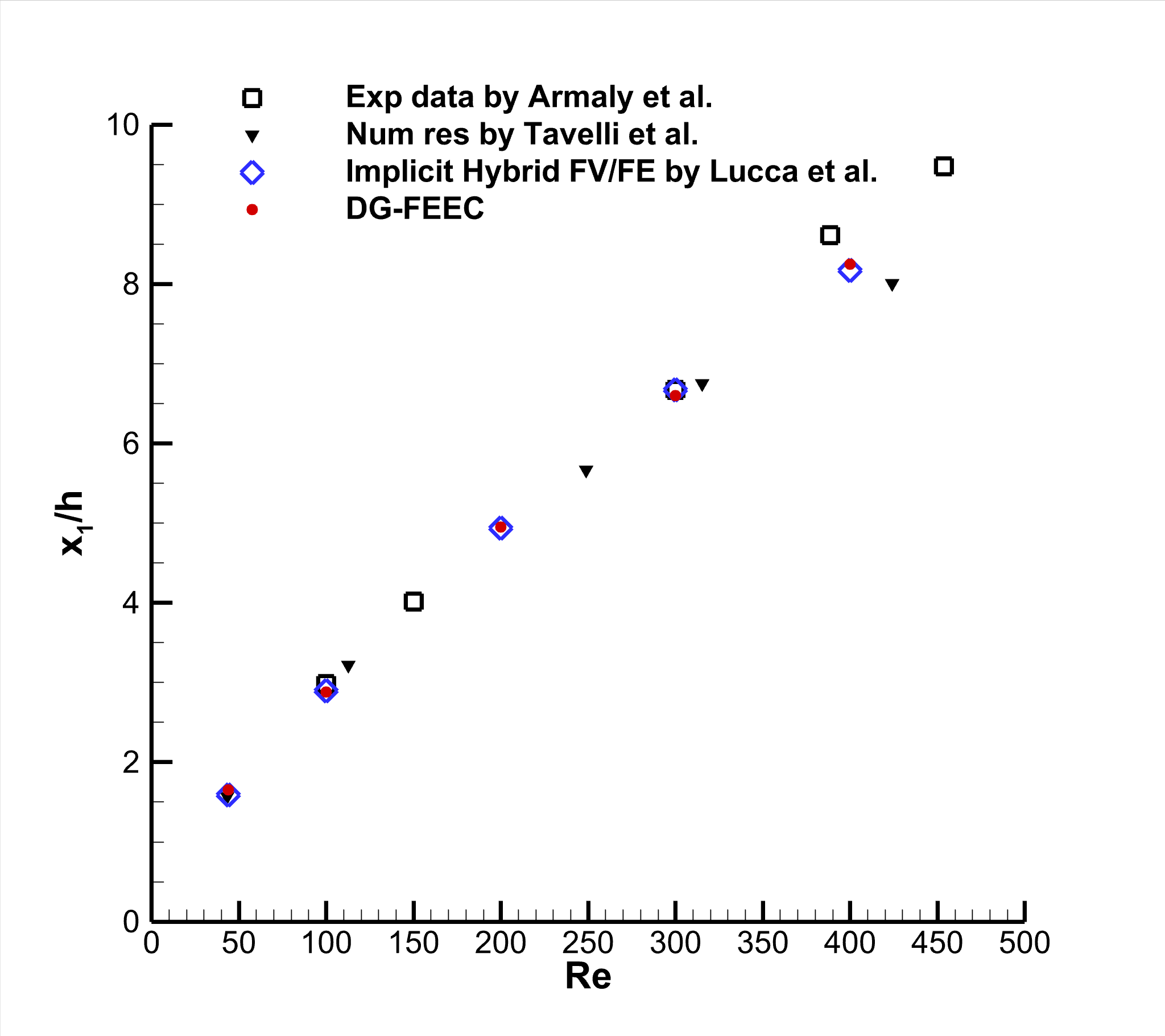}
	\caption{Streamlines and horizontal component of the velocity around the step for the backward-facing step test for $\mathrm{Re}= 400$ (left) and normalized recirculation point versus Reynolds number, compared with the experimental results from \cite{armaly1983} and the numerical ones from \cite{TD14,HybridImplicit} (right).}
	\label{fig:BFS}
\end{figure}

We consider now the backward-facing step problem originally investigated experimentally by Armaly et al. \cite{armaly1983}. From the numerical point of view, we follow the set up reported by Lucca et al. \cite{HybridImplicit}. In particular, the computational domain is $\Omega=[-L,0]\times[0,h] \cup [0,29.1]\times[-0.097,h]$ with $L = 1.94$ and $h = 0.103$. At the inlet we impose the Poiseuille velocity $\overline{\vel} = (\overline{u}_1,0)$ with
\begin{equation*}
	\overline{u}_1(x,y) = \frac{\Delta\pres}{2 L \mu}y(h-y),
\end{equation*} 
where $\Delta\pres = -3.060845359$. At the outlet we impose a constant pressure, while at all the other boundaries we impose no-slip wall boundary conditions. For this test, the Reynolds number $\mathrm{Re}$ is defined as 
\begin{equation*}
	\mathrm{Re} \doteq \frac{2hU}{\nu},
\end{equation*}
with $U$ being the average inlet velocity, i.e. 
\begin{equation*}
	U \doteq \frac{1}{h}\int_0^h \overline{u}_1\,dy.
\end{equation*}
We run the simulation until $\tEnd  = 80.0$ for $\mathrm{Re} = 44, 100, 200, 300, 400$. For this test, we set $C_{\mathrm{CFL}} = 0.1$. The results obtained with $r=2$ and $h_{\max} = 0.04$ are shown in Figure~\ref{fig:BFS}, in which we perform a qualitative comparison with the experimental results obtained by Armaly et al. \cite{armaly1983} and the numerical ones computed by Tavelli and Dumbser \cite{TD14} and by Lucca et al. \cite{HybridImplicit}.

\subsection{Viscous flow around a cylinder \textcolor{black}{at low Mach}}
\begin{figure}
	\centering
	\includegraphics[trim = {5 5 5 5}, clip, width = 0.55\textwidth]{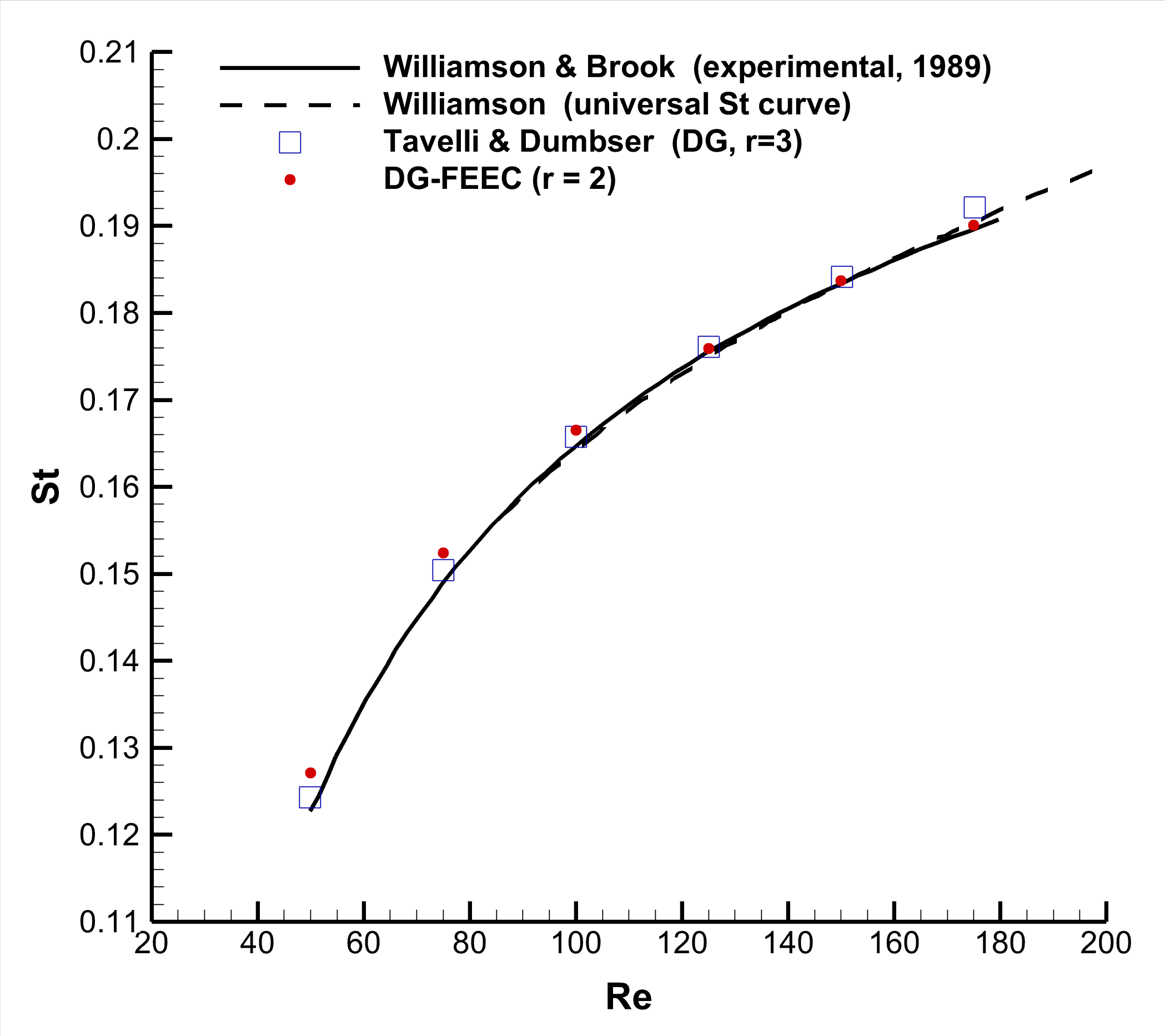}
	\caption{Strouhal number as a function of the Reynolds number.}
	\label{fig:ReVsSt}
\end{figure}

\begin{figure}
	\centering
	\includegraphics[trim = {0 0 0 1100}, clip, width = 0.75\textwidth]{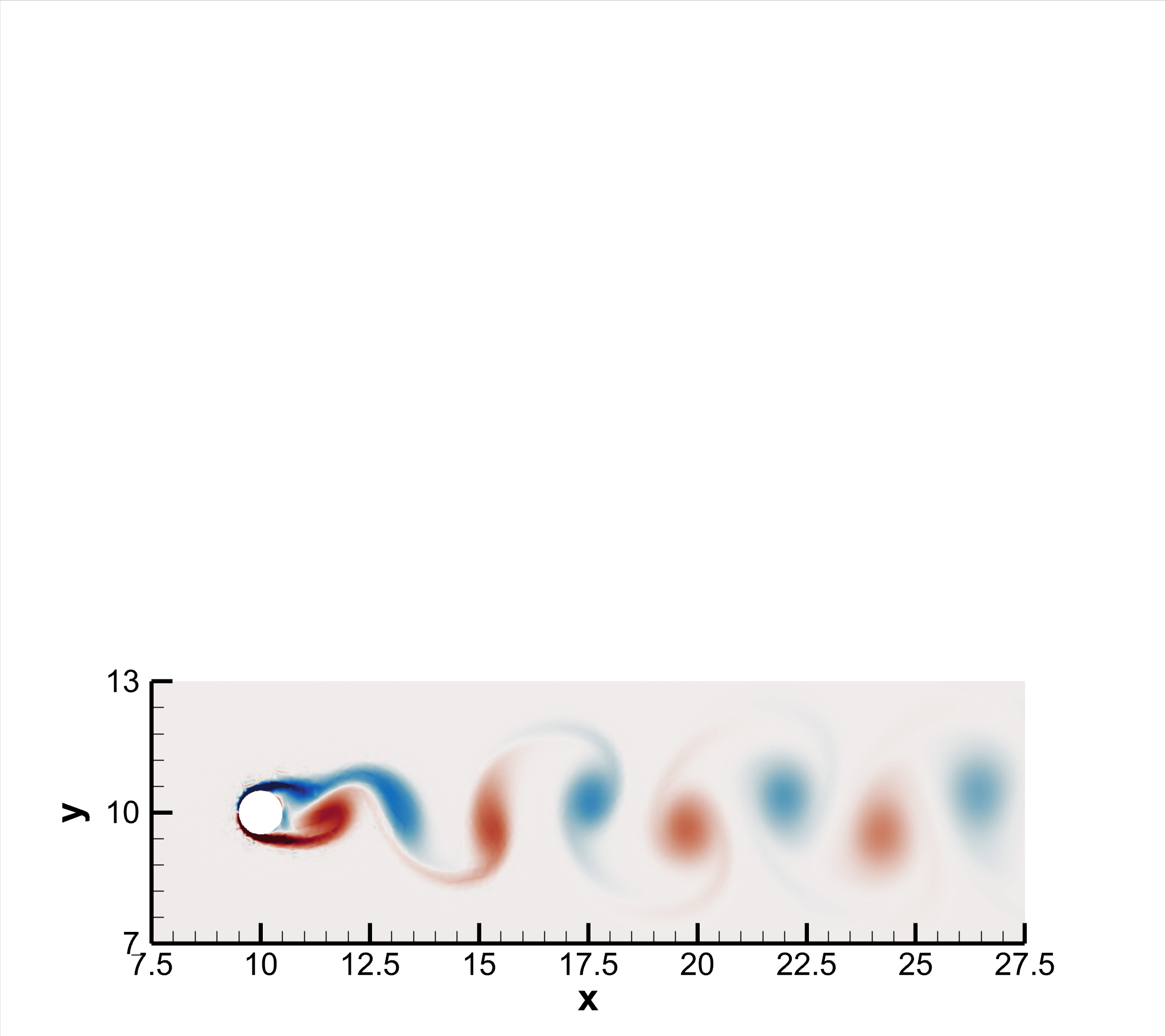}
	\caption{Detail of the vorticity around the cylinder for $\mathrm{Re} = 185$ at $t = 100$.}
	\label{fig:vonKarmann}
\end{figure}

We consider now the case of a viscous flow around a cylinder \cite{WillBrown98,HybridALE,HybridImplicit}. The computational domain is a rectangle with vertices $(0,0)$, $(50,0)$, $(50,20)$, $(0,20)$ minus a circle centered in $(10,10)$ with radius $0.5$.
The domain is discretized with a mesh made by 9168 elements with $\mathcal{P}_3$ curved boundaries around the cylinder. The following boundary conditions are imposed:
\begin{itemize}
	\item Inflow with $\overline{\vel} = (1,0)$ at the left boundary;
	\item Outflow with $\overline{p} = 0$ at all the other boundaries of the rectangle;
	\item No-slip wall on the boundary of the cylinder.
\end{itemize}
We compute the shedding frequency $f$ of the vorticity evaluated at the point $P = (15,10)$. In Figure~\ref{fig:ReVsSt} we plot the computed Strouhal number $\mathrm{St}$ (which for this test coincides with $f$) as a function of the Reynolds number $\mathrm{Re} = 1/\mu$. We compare our results with those obtained using the semi-implicit DG scheme proposed by Tavelli and Dumbser \cite{TD14}, the experimental data of Williamson and Brown \cite{WillBrown98} and the so-called universal Strouhal curve. \textcolor{black}{We remark that for this test we are using less elements than \cite{TD14} and \cite{HybridImplicit}, but the computed solution still agrees well with the reference solutions.} The vorticity at time $t = 100$ with $\mathrm{Re} = 185$ is shown in Figure~\ref{fig:vonKarmann}. 

\section{Conclusions}\label{sec:conclusions}
\color{black}
In this paper, we have introduced a novel semi-implicit method for weakly compressible flows based on compatible finite elements. This method achieves arbitrary high order in space and ensures exact mass conservation at the discrete level. Our proposed semi-implicit scheme leverages an operator splitting technique, as discussed in previous works \cite{Casulli1990,CasulliCheng1992,PM05,TV12,DC16}. The nonlinear convective terms are discretized using an explicit discontinuous Galerkin method, while all other terms are handled implicitly. Notably, each iteration involves solving only symmetric positive definite linear systems, thanks to the hybridization technique. When the Mach number approaches zero and density remains constant, our method tends to an exactly divergence-free scheme for the incompressible Navier-Stokes equations. The asymptotic-preserving property has also been verified numerically, where we find quadratic convergence of the density and the divergence errors in terms of the Mach number, as expected. Additionally, we have incorporated an \textit{a posteriori} limiter via artificial viscosity based on the MOOD approach. Finally, we validated the new scheme against a set of classical benchmark problems for both compressible and incompressible flows.

\color{black}
Due to the employed splitting approach, the numerical schemes proposed in this paper are so far limited to first order of accuracy in time. However, higher order time accuracy can be easily achieved, for example, at the aid of IMEX Runge-Kutta time integrators, see e.g. \cite{PareschiRusso2000,DLDV18,BDLTV2020,Thomann2020,Thomann2020b,Thomann2022,BDLTV2020}.  

\color{black}
Looking ahead, we plan to extend our scheme to viscous compressible flows by incorporating a discretization of the full Navier-Stokes tensor via the MCS method (e.g., \cite{GoLeSc20, GoLeSc20b}). Additionally, we aim to design a scheme capable of solving all Mach number flows, similar to the approaches in \cite{Hybrid2, TD17,  Thomann2022}. Another promising direction is extending our method to magnetohydrodynamics (MHD) by adding a conforming discretization of the magnetic field, as demonstrated in works such as \cite{HuMaXu17, HiLiMaZh18, GaGB22}.

%
\color{black}
\section*{Acknowledgements}
This work was financially supported by the Italian Ministry of Education, University 
and Research (MIUR) in the framework of the PRIN 2022 project \textit{High order structure-preserving semi-implicit schemes for hyperbolic equations} and via the Departments of Excellence  Initiative 2018--2027 attributed to DICAM of the University of Trento (grant L. 232/2016). The authors are member of the GNCS-INdAM (Istituto Nazionale di Alta Matematica) group. M.D. was also co-funded by the European Union NextGenerationEU (PNRR, Spoke 7 CN HPC). Views and
opinions expressed are however those of the author(s) only and do not necessarily reflect those of the European Union or the European Research Council. Neither the European Union nor the granting authority can be
held responsible for them. E.Z. is grateful to J. Sch\"{o}berl for the help with $\mathtt{NGSolve}$ and the stimulating discussions.

%
%
%
%
\bibliographystyle{elsarticle-num}
\bibliography{./biblio}


\end{document}